\documentclass[12pt]{article}
\usepackage{amsmath}
\usepackage{latexsym}
\usepackage{amssymb}
\usepackage{a4wide}
\newtheorem{thm}{Theorem}[section]
\newtheorem{la}[thm]{Lemma}
\newtheorem{Defn}[thm]{Definition}
\newtheorem{Remark}[thm]{Remark}
\newtheorem{Note}[thm]{Note}
\newtheorem{prop}[thm]{Proposition}

\newtheorem{cor}[thm]{Corollary}
\newtheorem{Example}[thm]{Example}
\newtheorem{Examples}[thm]{Examples}
\newtheorem{Problems}[thm]{Problems}

\newtheorem{Problem}[thm]{Problem}
\newtheorem{Number}[thm]{\!\!}
\newenvironment{defn}{\begin{Defn}\rm}{\end{Defn}}

\newenvironment{rem}{\begin{Remark}\rm}{\end{Remark}}
\newenvironment{numba}{\begin{Number}\rm}{\end{Number}}
\newenvironment{proof}{{\noindent\bf Proof.}}%
                  {\nopagebreak\hspace*{\fill}$\Box$\medskip\medskip\par}   
\newcommand{\Punkt}{\nopagebreak\hspace*{\fill}$\Box$}
\newcommand{\wb}{\overline}

\newcommand{\wt}{\widetilde}

\newcommand{\n}{\rm}

\newcommand{\ident}{\equiv}
\newcommand{\mto}{\mapsto}

\newcommand{\ve}{\varepsilon}
\newcommand{\isom}{\cong}

\newcommand{\Ad}{\mbox{{\rm Ad}}}
\newcommand{\ad}{\mbox{{\rm ad}}}
\newcommand{\Log}{\mbox{{\rm Log}}}

\newcommand{\N}{{\mathbb N}}

\newcommand{\R}{{\mathbb R}}

\newcommand{\Q}{{\mathbb Q}}

\newcommand{\Z}{{\mathbb Z}}

\newcommand{\K}{{\mathbb K}}

\newcommand{\F}{{\mathbb F}}
\newcommand{\A}{{\mathbb A}}

\newcommand{\cm}{{\mathfrak m}}

\newcommand{\wh}{\widehat}

\newcommand{\sub}{\subseteq}

\newcommand{\GL}{\mbox{\rm GL}}

\newcommand{\gl}{\mbox{gl}}

\newcommand{\car}{\mbox{{\rm char}}}

\newcommand{\id}{\mbox{\n id}}

\newcommand{\cL}{{\cal L}}

\newcommand{\sbull}{{\scriptscriptstyle \bullet}}

\begin{document}
%
%
%$\;$\vspace*{-2 cm}\\
\begin{center}
{\Large \bf Every smooth {\boldmath $p$\/}-adic
Lie group admits a compatible analytic structure}\vspace{4.1 mm}\\
{\bf Helge Gl\"{o}ckner}\vspace{.7mm}
\end{center}
\noindent{\bf Abstract.\/}
We show that every finite-dimensional
$p$-adic Lie group of class $C^k$ (where $k\in \N\cup\{\infty\}$)
admits a $C^k$-compatible analytic Lie group
structure.
We also construct an exponential map
for every $k+1$ times strictly differentiable
($SC^{k+1}$)
ultrametric $p$-adic Banach-Lie group, which is an $SC^1$-diffeomorphism
and admits Taylor expansions of all finite orders
$\leq k$.\\[3mm]
{\footnotesize
{\bf AMS Subject Classification.}
Primary 22E20, 22E65. %gen pro and struct of other Lie groups;
% inf-dim Lie groups
Secondary
22A05, % = structure of general topological groups
22D05, % general properties of lcp groups
22E35\\[3mm] % Analysis on $p$-adic Lie groups
{\bf Keywords and Phrases.}
$p$-adic Lie group, Banach-Lie group,
pro-$p$-group, filtration,
Taylor expansion,
strict differentiability,
smoothness, analyticity,
one-parameter group, $p$-adic interpolation,
exponential map, logarithm, canonical
coordinates, Lazard's theorem}
\section*{Introduction}
\noindent
It is well-known that every
finite-dimensional real Lie group
of class $C^k$ (where $1\leq k\leq \infty$)
admits a $C^k$-compatible
analytic manifold structure making the group
operations analytic
(\cite{Bir}, \cite{Eis},
\cite{Pon}, \cite{Sch}, \cite{Seg}, \cite{Smi};
also \cite[\S4.4]{MaZ}),
and similar (slightly  weaker) results are valid
for real Banach-Lie groups~\cite{Mai}.
It therefore suffices for all
practical purposes to ignore the case
of finite order differentiability
and restrict attention to
smooth real Lie groups,
or to analytic real Lie groups
whenever this is profitable.\\[3mm]
In this paper, we prove analogous
results for $p$-adic Lie groups:\\[3mm]
{\bf Theorem~A.}
{\em Let $G$ be a finite-dimensional Lie group of class
$C^k_\K$
$($where $k\in \N\cup \{\infty\})$
over a valued field $\K$
which is a finite extension of $\Q_p$.
Then there exists a $\K$-analytic manifold
structure on~$G$ making it a $\K$-analytic Lie group,
and which is $C^k_\K$-compatible with the given
$C^k_\K$-manifold structure.}\\[3mm]
Here, $C^k_\K$-maps are understood in the sense
of \cite{BGN}, where a setting
of differential calculus over arbitrary non-discrete
topological fields~$\K$
is described.
We recall that mappings between
open subsets of finite-dimensional
$\Q_p$-vector spaces are $C^1_{\Q_p}$ if and only
if they are strictly differentiable
(where strict differentiability is a classical
concept, \cite{BVa}).
Thus, taking
$k=1$, Theorem~A subsumes
as a special case
that every strictly differentiable
finite-dimensional $p$-adic Lie group
admits a compatible analytic structure
making it a $p$-adic Lie group in the usual sense,
as considered in \cite{BLi}, \cite{Laz}, and \cite{Ser}.
We also obtain
valuable structural information
concerning ultrametric Banach-Lie groups
(Proposition~\ref{filtration}, Remark~\ref{remfilt},
Proposition~\ref{powers}, Corollary~\ref{crexp}):\\[3mm]
{\bf Theorem~B.}
{\em Let $k\in \N\cup\{\infty\}$
and $\K$ be a complete ultrametric field.
Assume that
$G$ is a Lie group of class $C^{k+1}_\K$
modelled on an ultrametric Banach space
over~$\K$,
or that $\K$ is locally compact
and $G$ a finite-dimensional $C^k_\K$-Lie group
over~$\K$.
Then $G$ is complete,
and there exists a diffeomorphism $\phi\!: U\to B_r(0)$
from an open subgroup $U\sub G$ onto an open ball
$B_r(0)\sub L(G)$,
with the following properties\/}:
\begin{itemize}
\item[(a)]
{\em The inverse images $U_s:=\phi^{-1}(\wb{B}_s(0))$
of closed balls for $s\in \;]0,r[$
are open, normal subgroups of $U$.
The family of open subgroups
$(U_s)_{s<r}$
is a basis of identity neighbourhoods
for~$G$.
If $\,G$ is $C^2_\K$, then $(U_s)_{s<r}$
defines a filtration for $U$.}
\item[(b)]
{\em If $\car(\K)=p>0$,
then $\Z\to U$, $n\mto x^n$
extends to a continuous homomorphism $\eta_x\!: \Z_p\to U$, for
each $x\in U$.}
\item[(c)]
{\em If $\car(\K)=0$
and $\K$ is a valued extension
field of $\Q_p$ for some~$p$,
then
$\Z\to U$, $n\mto x^n$
extends to a continuous homomorphism $\eta_x\!: \Z_p\to U$
which actually is $C^1_{\Q_p}$,
for each $x\in U$.
The map
$\log_G\!: U\to L(G)$, $\log_G(x):=\eta_x'(0)$
is an $SC^1_\K$-diffeomorphism\/}\footnote{A map
between open subsets of normed spaces is called
$SC^1_\K$ if it is strictly
differentiable at each point. Every $C^2_\K$-map is $SC^1_\K$,
and every $SC^1_\K$-map is $C^1_\K$ \cite{IMP}.}
{\em onto $B_r(0)\sub L(G)$,
such that
\[
\exp_G:=\log_G^{-1}\!: B_r(0)\to U
\]
is an exponential map\/}
({\em see Definition\/}~\ref{defnexp})
{\em and an $SC^1_\K$-diffeomorphism.
In the Banach case, let us assume now that $G$
is not only $C^{k+1}_\K$
but $C^{k+2}_\K$ $($or at least $k+1$ times
strictly differentiable$)$.
Then
each continuous homomorphism $\Z_p\to G$
is $C^k_{\Q_p}$,
and we can achieve that, in local coordinates, $\exp_G$
admits Taylor expansions of all finite orders $\leq k$.}
\end{itemize}
{\bf Classical construction of the analytic structure.}
The exponential map $\exp_G\!: L(G)\to G$
of
a finite-dimensional smooth Lie group $G$
over $\R$
can be obtained via $\exp_G(X):=\gamma_X(1)$,
where $\gamma_X$ is the uniquely
determined integral curve starting in~$1$,
to the left invariant vector field $X_\ell$
on~$G$ with $X_\ell(1)=X$.
Here $\gamma_X$ can also be described as the unique
$C^1$-homomorphism $\R\to G$
with $\gamma'_X(0)=X$.
Due to smooth dependence of solutions on
parameters,
$\exp_G$ is smooth and
induces a local diffeomorphism at~$0$,
giving rise to a chart for $G$
(``coordinates of the first kind'').
It turns out that the group operations
are analytic with respect to this chart,
because solutions to linear differential
equations with analytic coefficients are analytic,
and depend analytically on initial conditions and
parameters. Using translations, the analytic structure
can then be spread over all of~$G$.\\[3mm]
{\bf Difficulties in the ultrametric case.}
It is impossible to adapt the classical construction
just described to the ultrametric case.
For example,
integral curves stop to be useful,
since neither can we assume their existence
(nor local uniqueness),
nor do they need to define
homomorphisms (nor local homomorphisms)
if they exist.
For instance,
consider a smooth injective map $\gamma\!: \Z_p\to\Q_p$
such that $\gamma(0)=0$ and $\gamma '\ident 0$
\cite[Exercise\,29.G]{Sch}.\footnote{Typical examples
of such maps send a $p$-adic integer
$\sum_{j=0}^\infty a_jp^j$ (where $a_j\in \{0,\ldots,p-1\}$)
to $\sum_{j=0}^\infty a_jp^{k_j}$,
where $0<k_1<k_2<\cdots$ is an ascending
sequence of integers going
to $\infty$ sufficiently fast.}
Any such $\gamma$ provides an example of a smooth integral curve for
the invariant vector field $X_\ell\ident 0$
on $(\Q_p,+)$,
such that $\gamma$
does not coincide with a homomorphism on any
zero-neighbourhood in $\Z_p$,
and fails to be analytic on any zero-neighbourhood.
Furthermore, the germ of $\gamma$ is not determined
by the initial value problem, because
also $\eta\ident 0$ is a solution.\footnote{Actually, the germs around~$0$
of smooth maps $\zeta\!:\Z_p\to \Q_p$
with $\zeta(0)=0$ and $\zeta'\ident 0$ form a
$\Q_p$-vector space
of dimension $2^{\aleph_0}$.}\\[3mm]
{\bf Strategy of proof.}
In view of the difficulties just described,
it is clear that a different strategy of proof
is needed to construct an analytic structure
on a finite-dimensional $p$-adic $C^k$-Lie group~$G$,
which strictly avoids
any recourse to differential equations.
Here, the crucial idea is to use
Lazard's characterization
of analytic $p$-adic Lie groups:\\[2.5mm]
{\em A topological group $G$ can be given a $($necessarily
unique$)$ finite-dimensional
$p$-adic analytic Lie group structure if and only if
$G$ has an open, compact subgroup~$U$ such that\/}:
\begin{itemize}
\item[{\bf L1}]
{\em $U$ is a pro-$p$-group};
\item[{\bf L2}]
{\em $U$ is finitely generated
topologically, i.e., $U=\wb{\langle F\rangle}$
for a finite subset $F\sub U$; and\/}:
\item[{\bf L3}]
{\em The set of $p^2$-th powers of elements of $U$
contains the commutator subgroup of~$U$,
$[U,U]\sub \{x^{p^2}\!: x\in U\}$.}
\end{itemize}
(Cited from \cite[p.\,157]{Ser};
cf.\ \cite[A1, Thm.\,(1.9)]{Laz}).
Hence, to construct a compatible
analytic structure on a finite-dimensional
$p$-adic $C^k$-Lie group $G$,
we need to master
two tasks:
\begin{itemize}
\item[(a)]
Show that $G$ satisfies conditions {\bf L1--L3}.
\item[(b)]
Show that Lazard's analytic Lie group structure
is $C^k$-compatible with the given $C^k$-manifold
structure.
\end{itemize}
Here, an open subgroup $U$
of $G$
satisfying {\bf L1} and {\bf L3}
can be constructed as in the analytic case,
based on Inverse
Function Theorems for $C^k$-maps (see \cite{IMP})
and Taylor expansions.
In the Banach case,
similar arguments also provide the
filtration $(U_s)_{s<r}$
described in Theorem~B\,(a),
and show that $\{x^{p^n}\!: x\in U_s\}\sub U_{p^{-n}s}$
for each $n\in \N$
in the situation of Theorem~B\,(b)
and (c), entailing that the homomorphism
$\Z\to U$, $z\mto x^z$
is continuous with respect to the $p$-adic topology
and hence extends to a continuous
homomorphism $\eta_x\!: \Z_p\to G$
(Section~2). We now assume that $\K$ is an extension field
of $\Q_p$.
To establish {\bf L2},
the natural strategy is to try and introduce
coordinates of the second kind on $U$.
So far, we only know that
{\em continuous\/} $p$-adic one-parameter
groups exist.
In a first step (Section~3),
we show that the one-parameter groups
$\eta_x$ are $C^1_{\Q_p}$ in fact,
which enables us to define $\log_G$
(as in Theorem~B\,(c)),
deduce its strict differentiability
by iterating first order
Taylor expansions of the $p$-th power map,
and to define $\exp_G:=\log_G^{-1}$.
While it is
not hard to see that every $C^1$-homomorphism
between smooth Lie groups modelled
on real locally convex spaces is automatically
smooth~\cite{DIR}, such a result is
not available here.\footnote{In the $p$-adic case,
one needs to control
iterated difference quotient maps
instead of mere higher differentials,
making it impossible to adapt the proof.}
It is therefore necessary
to prove by hand that every one-parameter group
$\Z_p\to G$
is~$C^k_{\Q_p}$.
This is the most difficult
problem arising in our construction.
We solve it in Section~4,
where we express
$\exp_G$ as a limit of iterates
of the inverse of the $p$-th power map~$\tau_p$:
\[
\exp_G(x)=\lim_{n\to\infty} \tau_p^{-n}(p^nx)
\]
and deduce that (in local coordinates)
$\exp_G$ admits Taylor expansions
of all finite orders $\leq k$
from the fact that so does~$\tau_p^{-1}$.
To this end, we substitute the $k$-th order
Taylor expansions of each $\tau_p^{-1}$ in the $n$-th
iterate $x\mto \tau_p^{-n}(p^nx)$
into each other, multiply out,
and then prove by tracking each individual
term of the resulting sum that
the homogeneous terms of each given order
converge, and that the
polynomials resulting as the
limits actually provide a Taylor expansion for
$\exp_G$ with a remainder term vanishing
of the required order.
Now $\exp_G$ possessing a $k$-th
order Taylor expansion,
so do the one-parameter groups
$\Z_p\to G$, $z\mto \exp_G(zx)$,
whence they are $C^k_{\Q_p}$-maps
by Schikhof's Converse to Taylor's
Theorem for Curves
(\cite[Thm.\,10.7]{Sh1}, \cite[Thm.\,83.5]{Sh2}).\\[3mm]
In the case of a finite-dimensional $p$-adic $C^k_{\Q_p}$-Lie group,
it is now straightforward to verify
condition {\bf L2},
whence a $p$-adic analytic Lie group
structure on $G$ exists.
It only remains to show its $C^k_{\Q_p}$-compatibility with
the given structure.
After adapting the
Trotter product formula
and other classical ideas from
the analytic case to $p$-adic
$C^k$-Lie groups, we prove compatibility
by showing the existence of coordinates of the
second kind which work simultaneously
for both the analytic and the given
$C^k$-Lie group structure (Section~5).\\[3mm]
In an appendix, we prove various lemmas
compiled in a preparatory Section~1.\\[3mm]
{\bf Directions for further research.}
Beyond Theorem~A,
one would expect that every
$C^{k+2}_\K$-Lie group~$G$ modelled
on an ultrametric Banach space over
a valued extension field $\K$ of $\Q_p$
admits a $C^k_\K$-compatible $\K$-analytic
Lie group structure.
Two steps are still missing to achieve
this goal:
1. Show that $\exp_G$ is not only an $SC^1_\K$-diffeomorphism
admitting finite order Taylor expansions,
but actually a $C^k_\K$-diffeomorphism.
2. Show that $\exp_G$ induces an isomorphism of groups
from some ball in $L(G)$, equipped
with the Baker-Campbell-Hausdorff multiplication,
onto a subgroup of~$G$.
Here, the quite complicated notion of higher
order differentiability makes it technically
difficult to perform Step~1.
Possibly, Step~2 might be based
on an adaptation of the
construction
of the analytic structure
on uniformly powerful pro-$p$-groups
via coordinates of the first kind,
as described in \cite[\S9.4]{Dix}.
Recall that
the finite extension fields of $\Q_p$
occurring in Theorem\,A
are precisely
the local fields\footnote{That is to say,
the totally disconnected, locally compact,
non-discrete topological fields.}
of characteristic~$0$ \cite{Wei}.
From this perspective, it is very natural to ask
whether Theorem\,A remains valid
for finite-dimensional
Lie groups over local fields of
positive characteristic.
However, due to the lack of natural,
preferred coordinate systems (and the impossibility
of an analogue of Lazard's Theorem)
in positive characteristic, it is not at all clear
how one might approach the problem in this case.
\section{Basic definitions and facts}
In this section,
we briefly recall some notations
and basic facts from
the papers \cite{BGN} and
\cite{IMP},
which are our basic references for
differential calculus over
non-discrete topological fields
(see also \cite{Ber}
and \cite{CON}).
Then, various (new) lemmas are formulated,
most of which are based on Taylor expansions
and generalize classical facts
from multivariable calculus over~$\R$.
We recommend to take
these lemmas on faith at this point;
if desired, the proofs can be looked up later
in Appendix~A.
\begin{defn}\label{defncek}
Let $\K$ be a non-discrete (Hausdorff)
topological field, $E$ and $F$ be (Hausdorff)
topological $\K$-vector spaces,
and $f\!: U\to F$ be a map on an open subset
of~$E$. We call $f$ $C^0$ if it is continuous;
it is $C^1$ if it is $C^0$ and if there exists a
(necessarily unique) continuous map
$f^{[1]}\!: U^{[1]}\to F$
on $U^{[1]}:=\{(x,y,t)\in U\times E\times \K\!:
x+ty\in U\}$ such that $f^{[1]}(x,y,t)=t^{-1}(f(x+ty)-f(x))$
if $t\not=0$.
Inductively, $f$ is called $C^k$ for $k\in \N$
if $f$ is $C^1$ and $f^{[1]}$ is $C^{k-1}$;
we then let $f^{[k]}:=(f^{[1]})^{[k-1]}$.
The map $f$ is $C^\infty$ or smooth if it is $C^k$ for all
$k\in \N$. We also write $C^k_\K$
instead of $C^k$ ($k\in \N_0\cup\{\infty\}$)
to emphasize the ground field.
\end{defn}
If $E=\K$ here, we call $f\!: U\to F$ a {\em curve\/}.
A curve $f\!: U\to F$ is $C^1$ if and only if there
exists a continuous map
$f^{<1>}\!: U\times U\to F$ such that $f^{<1>}(x,y)=(x-y)^{-1}(f(x)-f(y))
=:f^{>1<}(x,y)$
for all $x,y\in U$ such that $x\not=y$ \cite[La.\,6.1]{BGN}.
\begin{numba}
It can be shown that compositions of composable
$C^k$-maps are $C^k$ \cite[Prop.\,4.5]{BGN}.
Furthermore,
for every $C^k$-map $f\!: E\supseteq U\to F$,
the iterated directional derivatives
\[
d^jf(x,y_1,\ldots,y_j):=(D_{y_1}\cdots D_{y_j}f)(x)\qquad
\mbox{for $x\in U$, $y_1,\ldots,y_j\in E$}
\]
exist for all $j\in \N$ such that $j\leq k$,
and define continuous maps $d^jf\!:U\times E^j\to F$
such that $d^jf(x,\sbull)\!: E^j\to F$
is symmetric $j$-linear, for each $x\in U$.
The associated homogeneous
polynomials are denoted
$\delta^j_xf\!: E\to F$,
$\delta^j_xf(v):=d^jf(x,v,\ldots,v)$.
We occasionally abbreviate $f'(x):=df(x,\sbull)\!: E\to F$.
Every $C^k$-map $f\!: U\to F$ as before, where $k\in \N$,
admits a $k$-th order Taylor expansion
\[
f(x+ty)-f(x)=\sum_{j=1}^kt^j\, a_j(x,y)+t^kR_k(x,y,t)\qquad
\mbox{for all $(x,y,t)\in U^{[1]}$,}
\]
for uniquely determined
maps $a_j\!: U\times E\to F$ of class $C^{k-j}$
which are $F$-valued forms of degree $j$
in the second argument, and where $R_k\!: U^{[1]}\to F$
is a continuous map such that $R_k(x,y,0)=0$
for all $(x,y)\in U\times E$
(see \cite{BGN}, Thm.\,5.1 and Thm.\,5.4,
where $R_k$
is denoted $R_{k+1}$).
If $f$ is $C^\ell$
actually with $\ell>k$, then $R_k$ is $C^{\ell-k}$.
Since
$\delta^j_xf(y)=j!\,a_j(x,y)$,
the Taylor expansion attains the familiar form
when $\car(\K)=0$.
\end{numba}
\begin{numba}
If $\K$ is a valued field here
(whose absolute value $|.|\!: \K\to [0,\infty[$
we shall always assume non-trivial),
$E$ a normed $\K$-vector space
and $F$ a polynormed $\K$-vector space
(viz.\ a topological $\K$-vector space~$F$
whose vector topology can be obtained
from a family of continuous seminorms $\|.\|_\gamma$
on~$F$), then it is also natural to consider
\[
R\!: U\times U\to F,\qquad R(x,y):=R_k(x,y-x,1)\,.
\]
We call both $R_k$ and $R$ the $k$-th order Taylor remainder;
no confusion is likely.
\end{numba}
\begin{numba}
Let $E$ and $F$ be normed vector spaces over
a valued field~$\K$, and
$f\!: U\to F$ be a map on an open subset
of $E$. Given $x\in U$,
$f$ is called {\em strictly differentiable at $x$\/}
if $\lim_{(y,z)\to (x,x)}\frac{\|f(z)-f(y)-f'(x).(z-y)\|}{\|z-y\|}=0$
for some (necessarily unique) continuous linear map
$f'(x)\!: E\to F$ (where $y\not=z$).
The map $f$ is {\em strictly differentiable\/} (or $SC^1$)
if it is strictly differentiable at each $x\in U$.
Then $f$ is $C^1$, and $f'(x)=df(x,\sbull)$.
Recursively, the map is called $k$ times strictly differentiable
(or $SC^k$) if it is $SC^1$ and $f^{[1]}$ is $SC^{k-1}$.
The following can be shown (see \cite{IMP}):
If $f$ is $SC^k$, then $f$ is $C^k$.
If $f$ is $C^{k+1}$, then $f$ is $SC^k$.
If $\K$ is locally compact and $E$ is finite-dimensional,
then $f$ is $C^k$ if and only if $f$ is $SC^k$.
\end{numba}
\begin{numba}
If $(E,\|.\|)$ is a normed space over a valued
field $(\K,|.|)$, $x\in E$ and $r>0$,
we write $B_r^E(x):=\{y\in E\!: \|y-x\|<r\}$
and $\wb{B}_r^E(x):=\{y\in E\!: \|y-x\|\leq r\}$,
or simply $B_r(x)$ and $\wb{B}_r(x)$, when $E$ is understood.
If the absolute value
$|.|$ on $\K$ satisfies the ultrametric inequality
$|x+y|\leq\max\{|x|,|y|\}$,
we shall call $(\K,|.|)$ an {\em ultrametric field}.
In this case, if also $\|.\|$
is ultrametric, then each ball
$B_r(0)$ or $\wb{B}_r(0)$ around $0$
is an open and closed additive subgroup of~$E$.
Then $B_r(x)=x+B_r(0)$ is a coset and hence
also open and closed, and $B_r(x)=B_r(y)$ for each $y\in B_r(x)$
(and likewise for $\wb{B}_r(x)$).

While in the case of a real normed space
every non-zero vector~$v$ can be
normalized,
in general we cannot find $r\in \K$
such that $\|r v\|=1$
when working over a valued field $(\K,|.|)$
whose value group $|\K^\times|$
is a proper subgroup of $]0,\infty[$.
As a substitute for normalization,
in many of the proofs we shall
fix some element $a\in \K^\times$
with $|a|<1$,
and then consider $a^{-k}v$
where $k\in \Z$ is chosen
such that $|a|^{k+1}\leq \|v\|<|a|^k$.
\end{numba}
The following lemmas
will be needed later
(for mappings into Banach spaces).
\begin{la}\label{taylor}
Let $E$ be a normed space
over a valued field~$\K$,
$F$ be a polynormed $\K$-vector space,
and $f\!: U\to F$ be a $C^2_\K$-map on an open subset
$U\sub E$. Then, for every $x_0\in U$
and continuous seminorm $\|.\|_\gamma$ on $F$,
there are $\delta>0$ and $C>0$ such that $B_{2\delta}(x_0)\sub U$
and
\[
\|f(x+y)-f(x)-f'(x).y\|_\gamma \leq C\cdot \|y\|^2\quad \mbox{for all
$\,x\in B_\delta(x_0)$ and $y\in B_\delta(0)$.}
\]
\end{la}
\begin{la}\label{taylor2}
Let $(\K,|.|)$ be a valued field,
$E,F$ be normed spaces over~$\K$,
$H$ be a polynormed $\K$-vector space,
$U\sub E$ and $V\sub F$
open zero-neighbourhoods, $\|.\|_\gamma\!: H\to [0,\infty[$
a continuous seminorm,
and $f\!: U\times V\to H$
be a mapping of class~$C^2_\K$ such that,
for certain continuous linear maps $\lambda\!: E\to H$
and $\mu\!: F\to H$, we have
\begin{eqnarray}
f(x,0) &=& f(0,0)+\lambda(x)
\quad
\mbox{for all $\,x\in U\,$ and}\label{cc1}\\
f(0,y) &=& f(0,0)+\mu(y)\quad
\mbox{for all $\,y\in V$.}\label{cc2}
\end{eqnarray}
Then
there exists $\delta>0$ and a constant $C>0$ such that
$B_\delta^E(0)\sub U$, $B_\delta^F(0)\sub V$ and
\[
\|f(x,y)-f(0,0)-\lambda(x)-\mu(y)\|_\gamma
\leq C\|x\|\,\|y\|\quad \mbox{for all $\,x\in B_\delta^E(0)$
and $\,y\in B_\delta^F(0)$.}
\]
\end{la}
\begin{la}\label{smallo}
Let $E$ be a normed vector space over
a valued field $\K$, $F$ be a polynormed
$\K$-vector space, $f\!: U\to F$ be a map on an
open subset of~$E$, and $k\in \N$.
If $f$ is of class
$C^{k+1}_\K$ or if $f$ is of class~$C^k_\K$,
$\K$ is complete and $E$ finite-dimensional,
then the map
\[
U\to \cL^k(E,F)\,,\quad
x\mto d^kf(x,\sbull)
\]
is continuous.
\end{la}
Here, we equip the $\K$-vector space $\cL^k(E,F)$
of continuous $k$-linear maps $\beta\!: E^k\to F$
with the vector topology determined by the family
of seminorms $\|.\|_\gamma\!: \cL^k(E,F)\to[0,\infty[$
defined via
$\|\beta\|_\gamma:=\inf\{C>0\!: (\forall x_1,\ldots,x_k\in E)\;
\|\beta(x_1,\ldots, x_k)\|_\gamma\leq C\|x_1\|\cdots\|x_k\|\}$,
where $\|.\|_\gamma\!: F\to [0,\infty[$ ranges through the
continuous seminorms on $F$
(the double use of the symbol $\|.\|_\gamma$
should not create confusion).
We set $\cL(E,F):=\cL^1(E,F)$
and $\cL(E):=\cL(E,E)$.
\begin{la}\label{ntayl}
Let $E$ be a normed vector space over
a valued field $\K$ of characteristic~$0$, $F$ be a polynormed
$\K$-vector space, $k\in \N$
and $f\!: U\to F$ be a map on an
open subset of~$E$ such that $f$ is
$C^{k+1}_\K$, or such that $f$ is of $C^k_\K$,
$\K$ locally compact, and $E$ finite-dimensional.
Let $R\!: U\times U\to F$,
$R(x,y):= f(y)-f(x)-\sum_{j=1}^k\frac{1}{j!}\,
\delta^j_xf(y-x)$
be the $k$-th order Taylor remainder.
Then, for any $z\in U$ and any continuous
seminorm $\|.\|_\gamma$ on $F$, we have
\begin{equation}\label{remaind}
\lim_{(x,y)\to (z,z)} \frac{\|R(x,y)\|_\gamma}{\|x-y\|^k}=0
\quad \mbox{$($where $x\not=y)$.}
\end{equation}
\end{la}
The following partial converse will be used essentially:
\begin{la}\label{alongcurv}
Let $E$ be a normed vector space over an ultrametric field~$\K$
of characteristic~$0$,
$F$ be an ultrametric Banach space over~$\K$,
$k\in \N$,
and $f\!: U\to F$ be a function on an open subset
$U\sub E$. We assume that $f$ can be written in the form
\[
f(y)=f(x)+\sum_{j=1}^k a_j(x,y-x)+R(x,y)\quad \mbox{for all $x,y\in U$,}
\]
where
$a_j\!: U\times E\to F$ is a continuous mapping
for $j\in \{1,\ldots, k\}$
such that $a_j(x,\sbull)\!: E\to F$
is a homogeneous polynomial of degree~$j$
for each $x\in U$,
$\;R\!: U\times U\to F$ is a continuous
mapping
such that $R(z,z)=0$, and {\rm (\ref{remaind})}
holds for each $z\in U$ and the norm
$\|.\|_\gamma:=\|.\|$ on~$F$. Then $f\circ \eta$
is of class $C^k_\K$, for each $C^k_\K$-curve
$\eta\!: \K\supseteq I \to U$. The maps $a_j$, $j=1,\ldots,k$,
are uniquely determined.
\end{la}
It is unknown whether
any $f$ as in Lemma~\ref{alongcurv}
actually is~$C^k_\K$, beyond the case $E=\K$.
\begin{la}\label{multilin}
Let $E,F$ be normed vector spaces over an ultrametric field $\K$,
with $F$ ultrametric. Let $n\in \N$
and $\beta\!: E^n\to F$ be a continuous
$n$-linear map.
Then, for all $\ve>0$ and
$u_i,v_i\in E$
such that $\|u_i\|\leq 1$, $\|v_i\|\leq 1$, and $\|u_i-v_i\|\leq \ve$
for $i=1,\ldots, n$,
we have
\[
\|\beta(u_1,\ldots, u_n)-\beta(v_1,\ldots,v_n)\|\leq \|\beta\| \ve\,.
\]
\end{la}
\begin{proof}
$\|\beta(u_1,\ldots, u_n)-\beta(v_1,\ldots,v_n)\|
=\left\|\sum_{k=1}^n
\beta(v_1,\ldots,v_{k-1},u_k-v_k,u_{k+1},\ldots,
u_n)\right\|$
$\leq \max_k
\|\beta(v_1,\ldots,v_{k-1},u_k-v_k,u_{k+1},\ldots,
u_n)\|\leq \|\beta\|\ve$.
\end{proof}
See \cite{BGN}
for the definition of $C^k$-manifolds
over a topological field~$\K$.
Given $k\in \N\cup\{\infty\}$,
a {\em Lie group of class $C^k$\/} modelled
on a topological $\K$-vector space~$E$ is a
group $G$, equipped with a $C^k$-manifold structure
modelled on~$E$
with respect to which the group
multiplication and inversion are $C^k$.
We abbreviate $L(G):=T_1(G)$;
if $k\geq 3$, then $L(G)$
has a topological
Lie algebra structure
with bracket $[X,Y]:=[X_\ell,Y_\ell](1)$
obtained from the Lie bracket
of the corresponding left invariant
vector fields on~$G$.
\section{Filtrations and continuous one-parameter
groups in ultrametric Banach-Lie groups}
In this section,
based on Taylor expansions and versions
of the Inverse Function Theorem,
we construct a basis for the filter of identity neighbourhoods
in
an ultrametric Banach-Lie group~$G$,
consisting of open subgroups of~$G$
which are diffeomorphic to balls in $L(G)$.
In the process,
we collect more and more information concerning
the properties of these subgroups
and the behaviour of various
mappings on them.
This enables us
to prove the existence of a large supply of
continuous $p$-adic one-parameter groups
(which actually cover the above neighbourhoods).
The later sections
also hinge on the present investigations,
which give us sufficient
control over all relevant maps, their differentials
and Taylor expansions to establish
differentiability properties for
one-parameter groups, the logarithm, and
the exponential map.\\[3mm]
In the following,
$(\K,|.|)$ is a complete ultrametric
field, with valuation ring $\A:=\{z\in \K\!: |z|\leq 1\}$,
the maximal ideal $\cm:=\{z\in \K\!: |z|<1\}$ of~$\A$,
and residue field $\F:=\A/\cm$. We let $p:=\car(\F)$
be the characteristic of~$\F$.
If $\K$ is locally compact, we let $q:=[\A:\cm]$
be the module of~$\K$, and define
$a:=\max\{|z|\!: z\in \cm\}$;
thus $|\K^\times|=a^\Z$.
Starting with Proposition~\ref{filtration}\,(h),
in the case where $\car(\K)=0$,
we shall assume for the rest of the proposition and all
later results based on the proposition
that the absolute value $|.|$ of $\K$ is non-trivial
on~$\Q$. Then $0<|p|<1$, and after replacing
$|.|$ with a suitable power we may assume without loss
of generality that $|p|=p^{-1}$,
entailing that the closure of
$\Q$
in $\K$ is $\Q_p$, equipped with its usual absolute value.
\begin{prop}\label{filtration}
Let $k\in \N\cup\{\infty\}$, and
$E$ be an ultrametric Banach space over $\K$.
If $\K$ is locally compact and $E$ is finite-dimensional,
define $k_+:=k$; otherwise, define $k_+:=k+1$.
Let $G$ be a Lie group of class
$C^{k_+}_\K$ modelled on~$E$,
and $\|.\|$ be an ultrametric norm on $L(G)\isom E$
defining its topology.
Then there exists a
$C^{k_+}_\K$-diffeomorphism $\phi\!: U\to B_r(0)$
from an open subgroup $U$ of~$G$
onto
the open ball $B_r(0)\sub L(G)$
for some $r\in \,]0,1[$,
such that $\phi(1)=0$ and $T_1\phi=\id_{L(G)}$.
Via
\[
\mu(x,y):=x*y:=\phi(\phi^{-1}(x)\phi^{-1}(y))\quad \mbox{for $\;x,y\in B_r(0)$}
\]
we obtain a group multiplication on $B_r(0)$
making $\phi$ an isomorphism of $C^{k_+}_\K$-Lie groups.
Let $\alpha\in \,]0,1[$.
Then, replacing $\|.\|$
with a positive multiple
and shrinking $r$ if necessary,
it can be achieved that the following holds:
\begin{itemize}
\item[\rm (a)]
For each $x\in B_r(0)$ and $s\in \,]0,r]$,
we have $x* B_s(0)=B_s(0)* x=x+B_s(0)$.
In particular, the left uniform structure
on $(B_r(0),*)$, its right uniform structure,
and the uniform structure on $(B_r(0),+)$
all coincide, and $(B_r(0),*)$ is a complete
topological group $($whence so is $G)$.
\item[\rm (b)]
$V_s:=B_s(0)$ is an open, normal subgroup
of $(B_r(0),*)=:V_r$, for each $s\in \,]0,r]$.
\item[\rm (c)]
For every $\ve>0$, there exists $\delta\in \;]0,r]$
such that $[V_s,V_s]\sub V_{\ve s}$ for all
$s\in \;]0,\delta]$.
If $k_+\geq 2$, then furthermore
$[V_t,V_s]\sub V_{\alpha ts}\sub V_{ts}$,
for all $\,s,t\in\,]0,r]$.
\item[\rm (d)]
The inversion map $\iota\!: V_r\to V_r$,
$x\mto x^{-1}$
is a surjective isometry.
For each $\ve>0$, there exists $\delta\in \;]0,r]$
such that $\|x^{-1}+x\|\leq \ve \|x\|$ for all $x\in V_\delta$.
If $k_+\geq 2$, then furthermore
$\|x^{-1}+x\|\leq \alpha \|x\|^2$ for all $\,x\in V_r$.
\item[\rm (e)]
Define $W_s:= \wb{B}_s(0)$
for $s\in \,]0,r[$.
Then analogues of {\rm (b)\/} and {\rm (c)\/}
are valid
for $W_s$ and $W_t$ in place of $V_s$ and $V_t$
$($provided we require furthermore
$s,t <r$ here$)$.
\item[\rm (f)]
$V_s/V_{\alpha s}$ is abelian for each $s\in \,]0,r]$,
and so is $W_s/V_s$ for each $s\in \,]0,r[$.
\item[\rm (g)]
For all $\ve>0$, there exists $\delta\in\;]0,r]$
such that $\|xy-x-y\|\leq \ve\,\max\{\|x\|,\|y\|\}$
for all $x,y\in B_\delta(0)$,
writing $xy:=x*y$ now.
Hence $xy-x-y\in B_{\ve s}(0)$,
for all $s\in \;]0,\delta]$ and $x,y\in B_s(0)$.
If $k_+\geq 2$, we can achieve that
$\|xy-x-y\|\leq \alpha \|x\|\,\|y\|$
for all $x,y\in V_r$,
whence
$\|xy-x-y\|<\alpha s^2$
and thus $xy\in x+y+B_{\alpha s^2}(0)$,
for all $s\in\;]0,r]$ and $x,y\in V_s$.
\item[\rm (h)]
If $\K$ is locally compact,
$L(G)$ is of finite, positive dimension~$d$,
and $\|.\|$ corresponds to the maximum norm on $\K^d$
with respect to a suitable choice of basis
of $L(G)$, then we may achieve that $r=a^{j_0}$
for some $j_0\in \N$ and
\[
V_{a^j}/V_{a^{j+1}}\isom (\F^d,+)\quad\mbox{for all $j\in \N$ such that
$j\geq j_0$.}
\]
In particular, $[V_{a^j}:V_{a^{j+1}}]=q^d$
for each $j\geq j_0$,
and $U\isom V_r$ is a pro-$p$-group.
\item[\rm (i)]
Let
$\tau_p\!: V_r\to V_r$, $\tau_p(x):=x^p$
be the $p$-th power map.
If $\car(\K)=p$ and $k_+\geq 2$, we have
$(V_s)^{\{p\}}:=\tau_p(V_s) \sub V_{\alpha s^2}$
for all $s\in \;]0,r]$.
If $\car(\K)=p$ and $k_+=1$,
we have $(V_s)^{\{p\}}\sub V_{\alpha s}$
for each $s\in \;]0,r]$,
and furthermore
for each $\ve>0$ there exists $\delta\in \;]0,r]$
such that $(V_s)^{\{p\}}\sub V_{\ve s}$
for all $s\in \;]0,\delta]$.
If $\car(\K)=0$, we have
$(V_s)^{\{p\}}=V_{p^{-1}s}=pV_s$
for every $s\in \,]0,r]$,
and $p^{-1} \tau_p\!: V_r\to V_r$
is a surjective isometry.
\item[\rm (j)]
$[V_s,V_s]\sub V_{p^{-2}s}$, for all $s\in \,]0,r]$.
\item[\rm (k)]
$\|x^n\|\leq \|x\|$, for every $x\in V_r$ and $n\in \N_0$.
\item[\rm (l)]
For all $\ve>0$, there exists $\delta\leq r$
such that $\|x^n-nx\|\leq \ve \|x\|$ for all
$x\in V_\delta$ and all $n\in \Z$.
If $k_+\geq 2$, then furthermore
$\|x^n-nx\|\leq \alpha \|x\|^2$,
for every $x\in V_r$ and $n\in \Z$.
\item[\rm (m)]
For each $x\in V_r$, the homomorphism
$\Z\to V_r$, $n\mto x^n$ extends to
a continuous homomorphism $\eta_x\!: \Z_p\to V_r$.
Given $x\in V_r$ and $z\in \Z_p$, we abbreviate
$x^z:=\eta_x(z)$.
\end{itemize}
\end{prop}
\begin{proof}
(a) and (b): Let $\phi\!: P\to V \sub L(G)$ be a chart of~$G$
about~$1$, such that $\phi(1)=0$ and $T_1\phi=\id_{L(G)}$.
There is an open, symmetric identity neighbourhood $Q\sub P$
such that $QQ\sub P$.
Then $W:=\phi(Q)\sub V$ is an open zero-neighbourhood
in~$L(G)$, and
\[
\mu\!: W\times W\to V,\quad \mu(x,y):=x*y:=\phi(\phi^{-1}(x)\phi^{-1}(y))
\]
and
$\iota\!: W\to W$, $\iota(x):=\phi((\phi^{-1}(x))^{-1})$
are smooth maps
such that $\mu(x,0)=\mu(0,x)=x$,
$d\mu((0,0),(u,v))=u+v$, and $d\iota(0,u)=-u$
for all $x\in W$, $u,v\in L(G)$.
For each $x\in W$, the map
$\mu_x:=\mu(x,\sbull)\!: W\to V$
is a $C^{k_+}$-diffeomorphism onto its
open image.
Since $T_0(\mu_0)=\id_{L(G)}$,
there is
$r\in \;]0,1[$
such that $B_r(0)\sub W$,
$\mu_x(B_r(0))=B_r(0)$
for each $x\in B_r(0)$,
and
\begin{equation}\label{gp1}
\mu_x(B_s(0))=x+B_s(0)\qquad\mbox{for all $x\in B_r(0)$ and $s\in \,]0,r]$;}
\end{equation}
see \cite{IMP}, Thm.~7.4 (a)$'$ and (b)$'$.
Thus $x*B_s(0)=x+B_s(0)$.
After shrinking $r$, also
\begin{equation}\label{gpb}
B_s(0)*x=x+B_s(0)\quad
\mbox{for all $\,x\in B_r(0)$ and $s\in \,]0,r]$,}
\end{equation}
by an analogous argument.
Since $T_0(\iota)=-\id_{L(G)}$,
shrinking~$r$ even more if necessary,
we can achieve that
furthermore $\iota$ is isometric and
\begin{equation}\label{gp2}
\iota(B_s(0))=B_s(0)\qquad
\mbox{for all $s\in \;]0,r]$,}
\end{equation}
by \cite{IMP}, Prop.\,7.1\,(a)$'$ and (b)$'$.
Because $x+B_s(0)=B_s(0)$
for $x\in B_s(0)$, we deduce from
(\ref{gp1}) and (\ref{gp2})
that $B_s(0)$ is closed under the operations
$\mu$ and $\iota$ (which correspond to
multiplication and inversion in~$G$),
entailing that $\phi^{-1}(B_s(0))$
is a subgroup of $G$,
whence $V_s:=(B_s(0),\mu|_{B_s(0)\times B_s(0)},\iota|_{B_s(0)})$
is a Lie group, for each $s\in \;]0,r]$.
In particular, $U:=\phi^{-1}(B_r(0))$
is a subgroup of~$G$.
Since $x*V_s=x+B_s(0)=V_s*x$ for all
$x\in V_r$ by (\ref{gp1}) and (\ref{gpb}),
the subgroup $V_s$ of $V_r$ is normal.
The assertion concerning uniform structures
is clear from (\ref{gp1}) and (\ref{gpb}).
Now completeness transfers
from the open (and hence closed)
subgroup $(B_r(0),+)$ of $(L(G),+)$
to $(B_r(0),*)$, $U$, and $G$.

(c) Assume $k_+\geq 2$ first.
By Lemma~\ref{taylor2},
there exists a constant $C>0$ and $\delta<r$ such that
$\|xyx^{-1}y^{-1}\|\leq C\, \|x\|\,\|y\|$.
For the norm
$\|.\|':=\sigma \|.\|$, where $\sigma\geq \max\{1, C/\alpha\}$,
the preceding inequality turns into
$\|xyx^{-1}y^{-1}\|'\leq \frac{C}{\sigma}\|x\|'\|y\|'\leq
\alpha\|x\|'\|y\|'$. Hence, after replacing $\|.\|$ with $\|.\|'$,
the second assertion of (c) (and hence also the first)
is satisfied.

Case $k_+=1$: Then
$f\!: V_r\times V_r\to V_r$,
$f(x,y):=xyx^{-1}y^{-1}$ is strictly differentiable,
with $f(0,0)=0$ and
$f'(0,0)=0$ (using the maximum norm on
$L(G)\times L(G)$).
Given $\ve>0$, we therefore find $\delta\in \;]0,r]$ such that
$\|f(x,y)\|=\|f(x,y)-f(0,0)-f'(0,0).(x,y)\|\leq \ve \|(x,y)\|$
for all $x,y\in V_\delta$.
Then $f(V_s\times V_s)\sub V_{\ve s}$
and thus $[V_s,V_s]\sub V_{\ve s}$, for all $s\in \;]0,\delta]$.

(d) By the proof of (a) and (b), $\iota\!: V_r\to V_r$
is a surjective isometry. {\em Case $k_+\geq 2$\/}:
Since $\iota'(0)=-\id_{L(G)}$,
we deduce from Lemma~\ref{taylor}
that there are $\rho \in \;]0,r]$
and $C>0$ such that
$\|\iota(x)+x\|\leq C\|x\|^2$ for all $x\in V_\rho$.
Replacing $\|.\|$ with a suitable multiple (as in the proof
of (c)), we may assume that $C\leq \alpha$,
whence the third (and hence also the second)
assertion holds.
{\em If $k_+=1$}, then $\iota$ is strictly differentiable
with $\iota(0)=0$ and $\iota'(0)=-\id_{L(G)}$,
whence, given $\ve$, there is $\delta\in \;]0,r]$
such that $\|\iota(x)+x\|\leq \ve \|x\|$
for all $x\in V_\delta$.

(e) Using that $W_s=\bigcap_{b>s}V_b$ for each $s\in \,]0,r[$,
(e) readily follows from (b) and (c).

(f) Let $\delta$ be as in (c),
applied with $\ve:=\alpha$.
Then
$[V_s,V_s]\sub V_{\alpha s}$
for all $s\in \,]0,\delta]$,
whence $V_s/V_{\alpha s}$ is abelian.
Since $[W_s,W_s]\sub W_{\alpha s}\sub V_s$,
also $W_s/V_s$ is abelian. Now replace $r$ with~$\delta$.

(g) Case $k_+\geq 2$:
Since $\mu(x,0)=x$ and $\mu(0,y)=y$
for all $x,y\in V_r$,
applying Lemma~\ref{taylor2} with $\lambda=\mu=\id_{L(G)}$
we find $\rho\in \;]0,r]$
and $C>0$ such that
$\|xy-x-y\|\leq C \|x\|\, \|y\|$
for all $x,y\in V_\rho$.
After replacing $\|.\|$
with a suitable positive multiple, we may assume
that $C\leq \alpha$ (cf.\ proof of~(c)).
Replacing now
$r$ with $\rho$, we have
$\|xy-x-y\|\leq \alpha \|x\|\,\|y\|$
for all $x,y\in V_r$, whence the second assertion
holds and hence also the first.
If $k_+=1$, then $\mu$ is strictly differentiable
with $\mu(0,0)=0$ and $d\mu((0,0),\,(x,y))=x+y$,
from which the assertion
readily follows.

(h) By (g),
after shrinking $r$, we may assume that
$\|xy-x-y\|<a\,\max\{\|x\|,\|y\|\}$
for all $x,y\in V_r$, and
$r=a^{j_0}$
for some $j_0\in \N$.
The subgroup
$V_{a^{j+1}}$ being normal in $V_{a^j}$ by (b),
for each $j\geq j_0$,
we deduce that
\[
xV_{a^{j+1}}\, y V_{a^{j+1}}
=xyV_{a^{j+1}}=
xy+V_{a^{j+1}}=x+y+B_{a^{j+1}}(0)\qquad\mbox{for all $x,y\in V_{a^j}$}
\]
from (a) and the fact the
$\|xy-x-y\|<a^{j+1}$
by the preceding. 
Hence the quotient $V_{a^j}/V_{a^{j+1}}$ of $(V_{a^j},*)$
does not only coincide with the quotient
$B_{a^j}(0)/B_{a^{j+1}}(0)$ of $(B_{a^j}(0),+)$
as a set (see (a)), but also as a group.
Here
$B_{a^j}(0)/B_{a^{j+1}}(0)\isom (\A/\cm)^d\isom \F^d$.
Since $(V_{p^{-j}})_{j\geq j_0}$ is a basis of identity neighbourhoods
consisting of open normal
subgroups of $V_r$ such that $V_r/V_{p^{-j}}$ is a finite $p$-group,
indeed $V_r$ is a pro-$p$-group.

(i) If $\car(\K)=p$,
we have $\tau_p'(0)=p\, \id_{L(G)}=0$.
Thus Lemma~\ref{taylor} provides $\delta< r$
and $C>0$
such that $\|\tau_p(y)\|\leq C\|y\|^2$
for all $y\in B_\delta(0)$, when $k_+\geq 2$.
As in the proof of (c),
we see that $C$ can be chosen $\leq \alpha$
after replacing $\|.\|$
with a suitable positive multiple. Replacing
$r$ with $\delta$, the assertion then holds.
Since $\tau_p(0)=0$ and $\tau_p'(0)=0$,
the assertion concerning the case
$\car(\K)=p$ and $k_+=1$ readily follows from
the strict differentiability of~$\tau_p$
(possibly after shrinking $r$).
Thus (i) holds in positive characteristic.

Now assume $\car(\K)=0$.
Since $(\tau_p)'(0)=p\, \id_{L(G)}$,
there exists $\rho\in \;]0,r]$
such that $p^{-1}\tau_p\!: V_\rho\to V_\rho$
is a surjective isometry (\cite{IMP},
Prop.\,7.1 (a)$'$ and (b)$'$). Now replace $r$ by~$\rho$.

(j) Replace $r$ with $\delta$ from (c),
applied with $\ve:=p^{-2}$.

(k) Let $x\in V_r$ and $n\in \N$.
If
$x\in V_s$ for some $s\in \;]0,r]$,
then also $x^n\in V_s$, since $V_s$ is a group.
Hence $\|x^n\|=\inf\{s\!: x^n\in V_s\}\leq \inf\{s\!: x\in V_s\}=\|x\|$.

(l) Case $k_+\geq 2$: It suffices to prove
the second assertion. We assume $n\in \N_0$ first
and proceed by induction.
The case $n=0$ is trivial. Now assume the assertion is correct
for some $n\in \N_0$.
Then, exploiting the ultrametric inequality, we obtain
\begin{eqnarray*}
\|x^{n+1}-(n+1)x\| &=& \|x^nx-x^n-x+x^n-nx\|\\
&\leq &\max\{\|x^nx-x^n-x\|,\|x^n-nx\|\}
\leq \alpha\|x\|^2\,,
\end{eqnarray*}
using that $\|x^nx-x^n-x\|\leq \alpha\|x^n\|\,\|x\|\leq
\alpha\|x\|^2$ by (g) and (k),
and $\|x^n-nx\|\leq \alpha\|x\|^2$
by the induction hypothesis.

To complete the proof of the second assertion,
observe that
\[
\|x^{-n}+nx\|=\|(x^{-1})^n-nx^{-1}+nx^{-1}+nx\|\leq
\max\{\alpha \|x^{-1}\|^2,|n|\cdot \|x^{-1}+x\|\}\leq
\alpha\|x\|^2
\]
for each $n\in \N$,
using that $\|x^{-1}\|=\|x\|$ and $\|x^{-1}+x\|\leq \alpha\|x\|^2$
by (d), and $|n|\leq 1$ since $(\K,|.|)$ is an
ultrametric field.

Case $k_+=1$: Given $\ve>0$,
choose $\delta$ with the properties
described in (d) and (g).
Let $n\in \N$ first;
the case $n=1$ is trivial.
For $x\in V_\delta$, we have $\|x^n-nx\|\leq \ve\|x\|$
by induction and $\|x^nx-x^n-x\|\leq \ve\,\max\{\|x^n\|,\|x\|\}=\ve\|x\|$,
by (g) and (k). Thus
$\|x^{n+1}-(n+1)x\|=\|x^nx-x^n-x+x^n-nx\|
\leq \max\{\|x^nx-x^n-x\|,\|x^n-nx\|\}\leq \ve\|x\|$
indeed. Then also $\|x^{-n}+nx\|\leq \max\{
\|(x^{-1})^n-nx^{-1}\|,\|nx^{-1}+nx\|\}
\leq \max\{\ve \|x^{-1}\|,\|x^{-1}+x\|\}\leq \ve\|x\|$,
by the choice of $\delta$.

(m) Let $x\in V_r$.
Given $\ve\in \;]0,r]$, there exists $N\in \N$ such that
$\alpha^Nr<\ve$ and $p^{-N}r<\ve$.
Then $x^n\in V_\ve$
for each $n\in \Z$ such that $|n|_p\leq p^{-N}$.
Indeed, any such $n$ has the form
$n=p^jm$, where $j\geq N$ and $m$ is either $0$ or coprime
to $p$. Then $x^{p^j}=\tau_p^j(x)\in V_{\alpha^jr}\sub V_\ve$
if $\car(\K)=p$,
and $x^{p^j}=\tau_p^j(x)\in V_{p^{-j}r}\sub V_\ve$
if $\car(\K)=0$,
by (i).
Hence, $V_\ve$ being a group,
we also have $x^n=(x^{p^j})^m\in V_\ve$, as asserted.
By the preceding, the homomorphism
\begin{equation}\label{masp}
\Z\to V_r,\quad
n\mto x^n
\end{equation}
is continuous at $0$ with respect to
the topology on~$\Z$ induced by $\Z_p$ and therefore
uniformly continuous, being a homomorphism.
Now $(V_r,*)$ being complete by (a), we see
that the homomorphism from (\ref{masp})
extends uniquely to a continuous homomorphism
$\eta_x\!: \Z_p\to V_r$.
\end{proof}
\begin{rem}\label{remfilt}
In the preceding situation,
define $w\!: V_r\to \;]0,\infty]$
via $w(0):=\infty$,
$w(x):=\log_r(\|x\|)$
for $x\in V_r\setminus\{0\}$.
If $k_+\geq 2$,
then $w$ defines a filtration on $V_r\isom U$
(in the sense of \cite{Ser}, Part~I, Chapter~II,
Defn.\,2.1), by Proposition~\ref{filtration}\,(e).
The subgroups of $V_r$ associated to the filtration
(see \cite{Ser}, Part~I, Chapter~II, \S2)
are $(V_r)_\lambda=W_{r^\lambda}$, for each $\lambda\in \,]0,\infty[$.
\end{rem}
\section{Differentiability of one-parameter groups}
In this section,
we show that the continuous one-parameter groups
$\eta_x$ constructed above
(by $p$-adic interpolation
of power maps) are actually of class $C^1_{\Q_p}$.
This allows us to define a
logarithm and an exponential map,
which we then show to be
strictly differentiable.
\begin{la}\label{higherdhom}
If $\phi \!: G\to H$
is a homomorphism between $C^k$-Lie groups
over a non-discrete topological field and $\phi|_U$
is $C^k$ on some open
identity neighbourhood $U\sub G$, then $\phi$ is $C^k$.
\end{la}
\begin{proof}
For each $x\in G$, we have
$\phi|_{xU}=\lambda^H_{\phi(x)}\circ \phi|_U\circ \lambda_{x^{-1}}^G|_{xU}^U$
using the indicated left translations on $G$ and $H$,
which are $C^k$-maps. Hence $\phi|_{xU}$ is $C^k$.
\end{proof}
\begin{la}\label{highhom}
Let $\,G$ be a Lie group of class $C^1_\K$
over a non-discrete topological field~$\K$,
and $\xi\!: P\to G$ a continuous homomorphism,
defined on an open subgroup $P\sub \K$.
Assume that
\[
(\phi\circ \xi)'(0)=\lim_{z\to 0} \frac{\phi(\xi(z))}{z}
\]
exists for some
chart $\phi\!: U\to U_1\sub L(G)$ of~$G$
about~$1$ such that $\phi(1)=0$.
Then $\xi$ is $C^1_\K$.
\end{la}
\begin{proof}
Let $W\sub U$ be a symmetric open identity neighbourhood
such that $WW\sub U$, and $W_1:= \phi(W)$.
Define
$\mu\!: W_1\times W_1\to U_1$, $\mu(x,y):=x*y:=\phi(\phi^{-1}(x)\phi^{-1}(y))$.
Then $B:=\xi^{-1}(W)$ is an open $0$-neighbourhood in~$\K$;
we let $A\sub \K$ be an open, symmetric
$0$-neighbourhood such that $A+A\sub B$.
Define $\zeta\!: B\to W_1$, $\zeta(z):=\phi(\xi(z))$.
By hypothesis,
$\zeta'(0)=\lim_{z\to 0} \frac{\zeta(z)}{z}$
exists, whence
\[
\ve\!: B\to L(G)\, , \quad  \ve(z):=\left\{
\begin{array}{cl}
\frac{\zeta(z)}{z} &\;\mbox{if $z\not=0$}\\
\zeta'(0)&\;\mbox{if $z=0$}
\end{array}
\right.
\]
is a continuous map. For any $x,y\in A$
such that $x\not=y$, we have
\begin{eqnarray*}
\zeta^{>1<}(x,y)&=&\frac{1}{y-x}(\zeta(y)-\zeta(x))=
\frac{1}{y-x}(\zeta(x)*\zeta(y-x)-\zeta(x)*0)\\
&=&\mu^{[1]}\left( (\zeta(x),0),\bigl(0,(y-x)^{-1}\zeta(y-x)\bigr),
y-x\right)\\
&=& \mu^{[1]}((\zeta(x),0),(0,\ve(y-x)),y-x)\,,
\end{eqnarray*}
where the final expression also makes sense
for $x=y$, and defines a continuous function
$(\zeta|_A)^{<1>}\!: A\times A\to L(G)$. Thus $\zeta|_A$ is $C^1_\K$,
whence so is $\xi|_A$ and hence $\xi$, by Lemma~\ref{higherdhom}.
\end{proof}
\begin{prop}\label{powers}
If $\car(\K)=0$
in the situation of Proposition~{\rm \ref{filtration}},
then, shrinking $r$ if necessary,
it can be achieved that,
in addition to {\rm (a)--(m)},
also the following holds:
\begin{itemize}
\item[\rm (n)]
For all $\ve>0$, there exists $\delta\leq r$ such that
$\|x^z-zx\|\leq \ve\|x\|$, for all
$x\in V_\delta$ and $z\in \Z_p$.
If $k_+\geq 2$, then furthermore
$\|x^z-zx\|\leq \alpha\|x\|^2$, for all $x\in V_r$ and $z\in \Z_p$.
\item[\rm (o)]
For each $x\in V_r$, the homomorphism $\eta_x$
is of class $C^1_{\Q_p}$, and $\eta_x'(0)\in V_r$.
\item[\rm (p)]
The mapping $V_r\to \cL(L(G))$, $x\mto \tau_p'(x)$
is continuous, with $\tau_p'(0)=p\,\id_{L(G)}$
and\linebreak
$\|\tau_p'(x)-p\,\id_{L(G)}\|<p^{-1}$
for each $x\in V_r$, whence
$p^{-1}\tau_p'(x)$ is a surjective
isometry.
\item[\rm (q)]
The map
$\log\!: V_r\to V_r$, $\log(x):=\eta_x'(0)$
is a surjective isometry
and an $SC^1_\K$-diffeomorphism such that $\log(0)=0$
and $\log'(0)=\id_{L(G)}$. We have $\log(x^z)=z\log(x)$
for all $x\in V_r$ and all $z\in \Z_p$.
\item[\rm (r)]
The map $\exp:=\log^{-1}\!: V_r\to V_r$
is a surjective isometry and an $SC^1_\K$-diffeomorphism
satisfying $\exp(0)=0$ and $\exp'(0)=\id_{L(G)}$.
For every $x\in V_r$, the map $\zeta\!:
\Z_p\to V_r$, $z\mto \exp(zx)$
is a homomorphism of class $C^1_{\Q_p}$
such that $\zeta'(0)=x$, and it is uniquely determined
by this property.
\end{itemize}
\end{prop}
\begin{proof}
(n) In view of the continuity of $\xi_x$ and scalar multiplication
$\Z_p\times L(G)\to L(G)$, assertion (n)
readily follows from Proposition~\ref{filtration}\,(l).

(o)
The limit
$\lim_{z\to 0}z^{-1}x^z$ (with $0\not=z\in \Z_p$)
will exist uniformly in $x\in V_r$, provided the following limit exists
uniformly in $x\in V_r$:
\begin{equation}\label{simplcation}
\lim_{n\to \infty}p^{-n}x^{p^n}\,.
\end{equation}
To see this, let $\ve>0$.
By (n), there exists $\delta\in \;]0,r]$
such that
$\|x^z-zx\|\leq \ve p^{-1}\|x\|\leq \ve \|x\|$ for all $x\in V_\delta$
and $z\in \Z_p$.
We may assume that $\delta=p^{-N}$ for some $N\in \N$.
Let $x\in V_r$ and $0\not=z\in \Z_p$
with $|z|_p\leq \delta$.
Then $|z|_p=p^{-n}$ with $n\geq N$
and thus $z=p^nm$, where $m:=p^{-n}z\in \Z_p$ has absolute value~$1$.
Since $\|x^{p^n}\|=p^{-n}\|x\|<|z|_p r\leq \delta$,
we have $\|(x^{p^n})^m-mx^{p^n}\|\leq \ve \|x^{p^n}\|=\ve p^{-n}\|x\|$.
Noting that $x^z=(x^{p^n})^m$ (which is clear
if $m\in \Z$ and follows for $m\in \Z_p$
by continuity), we see that
$z^{-1}x^z=p^{-n}x^{p^n}+z^{-1}(x^z-mx^{p^n})$,
where $\|z^{-1}(x^z-mx^{p^n})\|\leq |z|_p^{-1} \ve p^{-n}\|x\|\leq \ve$
and $n\geq N$.
Thus indeed $z^{-1}x^z$ will converge uniformly
if so does $p^{-n}x^{p^n}$.
Let us show now that $p^{-n}x^{p^n}$ converges uniformly.
Given $\ve$,
with $N$ as before we have
$\|x^{p^{n+1}}-px^{p^n}\|\leq \ve p^{-1}\|x^{p^n}\|$
for all $n\geq N$ and thus
\begin{equation}\label{telescope}
\|p^{-n-1}x^{p^{n+1}}-p^{-n}x^{p^n}\|\leq \ve p^n \|x^{p^n}\|\leq \ve\quad
\mbox{for all $x\in V_r$ and all $n\geq N$.}
\end{equation}
Writing $p^{-n_2}x^{p^{n_2}}-
p^{-n_1}x^{p^{n_1}}$ as a telescopic sum,
(\ref{telescope}) and the ultrametric inequality
yield
\[
\left\|
p^{-n_2}x^{p^{n_2}}-
p^{-n_1}x^{p^{n_1}}
\right\|
\leq
\max_{j=0,\ldots, n_2-n_1-1}
\left\|
p^{-n_1-j-1}x^{p^{n_1+j+1}}-
p^{-n_1-j}x^{p^{n_1+j}}
\right\| \leq \ve
\]
for all $x\in V_r$
and $n_1,n_2\geq N$ (where $n_2>n_1$, say),
whence
$p^{-n}x^{p^n}$ is uniformly Cauchy and hence
converges uniformly in~$x$.

By the preceding, the limit
$\eta_x'(0)=\lim_{z\to 0}z^{-1}x^z=
\lim_{n\to\infty}p^{-n}x^{p^n}$ exists
for each $x\in V_r$. Since
$\|p^{-n}x^{p^n}\|=\|x\|$ for all~$n$,
also $\|\eta_x'(0)\|\leq \|x\|<r$,
whence $\eta_x'(0)\in V_r$.

(p) The map $\tau_p$ is strictly differentiable
in both of the cases $k_+=1$ and $k_+\geq 2$,
whence $x\mto \tau_p'(x)$ is continuous
by \cite[La.\,3.2]{IMP}.
Since $\tau_p'(0)=p\,\id_{L(G)}$,
after shrinking $r$
we therefore have
$\|\tau_p'(x)-p\,\id_{L(G)}\|<p^{-1}$
for each $x\in V_r$ and thus
$\|p^{-1}\tau_p'(x)-\id_{L(G)}\|<1$,
whence
$p^{-1}\tau_p'(x)$ is a surjective
isometry by \cite{IMP},
La.\,7.2 and its proof.

(q) Let $x\in V_r$ be given; we want to show
that $\log$ is strictly differentiable at~$x$.
To this end, assume
$\ve>0$. 
The map
\[
h\!: V_r\times V_r\to L(G)\,,\quad h(y,u):=\tau_p(y+u)
\]
being $C^2_\K$ (if $k_+\geq 2$),
resp., $C^1_\K$ (if $k_+=1$),
applying \cite[La.\,3.5]{IMP}
(resp., \cite[La.\,4.5]{IMP})
to $h$
around $(0,0)$,
we find $\delta\in \;]0,r]$
such that
\begin{equation}\label{reuseeq}
\| \tau_p(w)-\tau_p(v)-\tau_p'(y).(w-v)\|\leq \frac{\ve}{p}\, \|w-v\|
\end{equation}
for all $y,v,w\in V_\delta$.
There exists $N\in \N$ such that $\tau_p^N(x)\in V_\delta$.
Applying \cite[La.\,3.5]{IMP} (resp.,
\cite[La.\,4.5]{IMP})
to $h$ around $(\tau_p^{N-1}(x),0)$,
we find $\delta_{N-1}\in \;]0,r]$
such that (\ref{reuseeq}) holds
for all $y,v,w\in B_{\delta_{N-1}}(\tau_p^{N-1}(x))$.
After shrinking $\delta_{N-1}$,
we can achieve that $\tau_p(B_{\delta_{N-1}}(\tau_p^{N-1}(x)))
\sub V_\delta$ here. Proceeding in this way,
we find $\delta_0,\delta_1,\ldots,\delta_{N-1}$
such that
\[
\tau_p(B_{\delta_n}(\tau_p^n(x)))\sub B_{\delta_{n+1}}(\tau^{n+1}_p(x))
\]
for all $n\in \{0,1,\ldots , N-2\}$, and such that (\ref{reuseeq})
holds for all $n\in \{0,1,\ldots, N-1\}$
and all $y,v,w\in B_{\delta_n}(\tau_p^n(x))$.

Let $v,w\in B_{\delta_0}(x)$.
Then $\tau_p^n(v)\in B_{\delta_n}(\tau_p^n(x))$ for all $n\in\{0,\ldots, N-1\}$
and thus furthermore $\tau_p^n(v)\in V_\delta$
for all $n\geq N$, and likewise for $\tau_p^n(x)$ and
$\tau_p^n(w)$.
By (\ref{reuseeq}), we have
\[
\tau_p(w)-\tau_p(v)-\tau_p'(x).(w-v)=:r_1,\quad
\quad\mbox{where $\|r_1\|\leq \frac{\ve}{p}\|w-v\|$.}
\]
Thus
$\|p^{-1}\tau_p(w)-p^{-1}\tau_p(v)-A_1.(w-v)\|
\leq \ve \,\|w-v\|$, where
$A_1:=p^{-1}\tau_p'(x)$.
More generally,
$r_n:=\tau_p^n(w)-\tau_p^n(v)-\tau_p'(\tau_p^{n-1}(x)).(\tau_p^{n-1}(w)-\tau_p^{n-1}(v))$
satisfies
\begin{equation}\label{nov}
\|r_n\|\leq \frac{\ve}{p}\|\tau_p^{n-1}(w)
-\tau_p^{n-1}(v)\|= \frac{\ve}{p^n} \|w-v\|
\end{equation}
for all $n\in \N$,
by (\ref{reuseeq}) and (i). Then
\begin{equation}\label{general}
\|p^{-n}\tau_p^n(w)-p^{-n}\tau_p^n(v)-A_n.(w-v)\|
\leq \ve \,\|w-v\|\quad \mbox{for all $\, n\in \N$,}
\end{equation}
where
$A_n :=p^{-n} \bigl(\tau_p'(\tau_p^{n-1}(x))\circ\tau_p'(\tau_p^{n-2}(x))\circ
\cdots\circ \tau_p'(x)\bigr)
\in \cL(L(G))$.
In fact,
\begin{eqnarray*}
\tau_p^n(w)-\tau_p^n(v)
&=&
\tau_p'(\tau_p^{n-1}(x)).(
\tau_p^{n-1}(w)-\tau_p^{n-1}(v))
+r_n\\
&=&
\tau_p'(\tau_p^{n-1}(x)).
\tau_p'(\tau_p^{n-2}(x)).(\tau_p^{n-2}(w)-\tau_p^{n-2}(v))+\tau_p'(\tau_p^{n-1}(x)).r_{n-1}+r_n\\
&=& p^nA_n(w-v)+\sum_{k=1}^{n-1}\bigl(
\tau_p'(\tau_p^{n-1}(x))\circ\cdots\circ \tau_p'(\tau_p^{n-k}(x))\bigr).r_{n-k}
+r_n
\end{eqnarray*}
shows that
$p^{-n}\tau_p^n(w)-p^{-n}\tau_p^n(v)-A_n.(w-v)$ equals
\[
\sum_{k=1}^{n-1}p^{-k}\bigl(
\tau_p'(\tau_p^{n-1}(x))\circ\cdots\circ \tau_p'(\tau_p^{n-k}(x))\bigr).p^{-(n-k)}r_{n-k}
+p^{-n} r_n\,,
\]
where each summand involving
$r_{n-k}$ has norm $\leq \ve \|y-x\|$,
because $\|r_{n-k}\|\leq \frac{\ve}{p^{n-k}}\|y-x\|$
and the derivatives
of $\tau_p$ have norm $\leq p^{-1}$, by (p).
Since also $p^{-n}r_n$ has norm $\leq \ve\|y-x\|$
by (\ref{nov}),
the ultrametric inequality shows that
(\ref{general}) holds.

We shall presently see that the sequence $(A_n)_{n\in\N}$
converges in $(\cL(L(G)),\|.\|)$.
Let $A\in \cL(L(G))$ denote the limit.
Recall from (\ref{simplcation})
that $p^{-n}\tau_p^n(w)\to \log(w)$
and $p^{-n}\tau_p^n(v)\to \log(v)$ as $n\to \infty$.
Letting pass $n\to\infty$
in (\ref{general}), we find that
\[
\|\log(w)-\log(v)-A.(w-v)\|\leq \ve\, \|w-v\|\,.
\]
As $v,w \in B_{\delta_0}(x)$ were arbitrary,
this means that $\log$ is indeed strictly differentiable
at~$x$, with $\log'(x)=A$.

\noindent
To see that the sequence $(A_n)_{n\in \N}$ converges,
let $\ve\in \;]0,1[$.
Part (p) provides $\delta\in \;]0,r]$
such that $\|p^{-1}\tau_p'(y)-\id\|< \ve$
for all $y\in V_\delta$. Let $N\in \N$ such that $p^{-N}\leq \delta$.
Then $\|p^{-1}\tau_p'(\tau_p^n(x))-\id\|<\ve$
for each $n\geq N$, whence
$p^{-1}\tau_p'(\tau_p^n(x))$
is contained in the ball
$B_\ve(\id)\sub \cL(L(G))$,
which is an open subgroup of $\GL(L(G))$.
As a consequence, for
any $n_1,n_2\geq N$,
with $n_2>n_1$, we have
$A_{n_2}\circ A_{n_1}^{-1}=
p^{-1}\tau_p'(\tau_p^{n_2-1}(x))\circ\cdots\circ
p^{-1}\tau_p'(\tau_p^{n_1}(x))\in B_\ve(\id)$.
Thus $(A_n)_{n\in \N}$ is a Cauchy sequence
in the Banach-Lie group $\GL(\cL(G))$
(which is complete by (a)),
and hence convergent.

If $x=0$ here, then $\tau_p^n(x)=0$ for each
$n$ and thus $A_n=\id$ for each $n$,
entailing that $\log'(0)=\lim_{n\to\infty}A_n=\id$.
Using \cite{IMP}, Prop.\,7.1\,(a)$'$, (b)$'$
and Thm.\,7.3,
we therefore find $\rho\in \;]0,r]$
such that $\log(V_\rho)=V_\rho$
and such that $\log\!: V_\rho\to V_\rho$ is a surjective
isometry and an $SC^1_\K$-diffeomorphism.
Now replace $r$ with $\rho$.

Since $\eta_0\ident 0$,
it is clear that $\log(0)=\eta_0'(0)=0$.
Finally, given $x\in V_r$,
there exists $y\in V_r$ such that $\eta_y'(0)=\log(y)=x$.
For $n,m\in \Z$,
we have $(y^n)^m=y^{nm}$,
whence $\eta_{y^n}(z)=\eta_y(nz)$
for all $z\in \Z_p$ by continuity
of $\eta_y$ and $\eta_{y^n}$,
using that $\Z$ is dense in $\Z_p$.
Thus $\log(y^n)=\eta_{y^n}'(0)=n\eta_y'(0)=n\log(y)$
for all $n\in \Z$
and hence $\log(y^z)=z\log(y)$
for all $z\in \Z_p$, by continuity.

(r) It is immediate from (q) that $\exp$ is a surjective isometry,
an $SC^1_\K$-diffeomorphism, $\exp(0)=0$, and $\exp'(0)=\id$.
Given $x\in V_r$ and $z\in \Z_p$,
we have $zx=z\log(\exp(x))=\log(\exp(x)^z)$
and thus $\exp(zx)=\exp(x)^z\!=\eta_{\exp x}(z)=:\zeta(z)$,
which is a $C^1_{\Q_p}\!$-homomorphism
$\Z_p\to V_r$ such that $\zeta'(0)=\log(\exp(x))=x$.
Suppose that also $\xi\!: \Z_p\to V_r$ is
a $C^1_{\Q_p}$-homomorphism such that
$\xi'(0)=x$. Set $y:=\xi(1)$.
Then $\xi(n)=y^n=\eta_y(n)$
for all $n\in \Z$ and thus $\xi=\eta_y$ by continuity,
entailing that $x=\xi'(0)=\eta_y'(0)=\log(y)$.
Hence $\exp(x)=y$ and $\zeta(z)=\exp(x)^z=y^z=\xi(z)$
for all $z\in \Z_p$, whence
$\xi=\zeta$.
\end{proof}
\section{Higher differentiability of one-parameter groups}
In this section, we
perform the most difficult step of
our construction:
We show that the exponential map $\exp\!: V_r\to V_r$
(from Proposition~\ref{powers}\,(r))
admits Taylor expansions of all finite orders $\leq k$.
As a consequence,
every $p$-adic
one-parameter subgroup of $G$ will be~$C^k$.\\[3mm]
We start with a simple observation:
\begin{la}
Let $X$, $Y$ be metric spaces
and $f_n\!: X\to Y$ homeomorphisms for $n\in \N$,
which converge uniformly to a surjective isometry
$f\!: X\to Y$.
Then $f_n^{-1}\to f^{-1}$ uniformly.
\end{la}
\begin{proof}
For each $y\in Y$, we have, using that $f$ is an isometry:
\[
d(f^{-1}(y),f_n^{-1}(y))=d(y,f(f_n^{-1}(y)))=
d(f_n(f_n^{-1}(y)),f(f_n^{-1}(y)))
\leq \sup_{x\in Y}d(f_n(x),f(x))\,.
\]
The assertion follows.
\end{proof}
Since $p^{-n}\tau_p^n(x)\to \log(x)$ uniformly
on $V_r$ by the proof of Proposition~\ref{powers}\,(o),
where $\log\!: V_r\to V_r$ is a surjective isometry, the
preceding lemma shows that
\begin{equation}\label{explimit}
\exp(x)=\lim_{n\to\infty} \tau_p^{-n}(p^nx)\quad\mbox{uniformly in $x\in V_r$,}
\end{equation}
whence $\exp(x)=\lim_{n\to\infty}g^n(p^nx)$
with $g:=(\tau_p)^{-1}\!: V_{p^{-1}r}\to V_r$,
and $g^n:=\overbrace{g\circ \cdots \circ g}^{n}\!: V_{p^{-n}r}\to V_r$.\\[3mm]
Having expressed $\exp$ in terms
of a limit of iterates of a given map,
the following proposition
(our most difficult technical result)
establishes Taylor expansions
(and hence higher differentiability properties)
for $\exp$.
\begin{prop}\label{difficult}
Let $(\K,|.|)$ be an ultrametric field
extending $(\Q_p,|.|_p)$.
Let $(E,\|.\|)$ be an ultrametric Banach space over~$\K$,
$r>0$, $k\in \N$,
and $g\!: B_r(0)\to B_{pr}(0)$ be a
diffeomorphism
of class $C^{k+1}_\K$ $($resp., of class $C^k_\K$
if $\K$ is locally compact and $E$ finite-dimensional$)$,
such that $g(0)=0$,
$g'(0)=p^{-1}\id_E$,
the map $B_r(0)\to B_r(0)$, $x\mto pg(x)$ is a surjective isometry,
and such that\vspace{-2mm}
\begin{equation}\label{leq1}
\|g'(x)\|=p \quad \mbox{for all $x\in B_r(0)$}
\end{equation}
$($where $B_r(0),B_{pr}(0)\sub E)$.
Suppose that the surjective isometries
\[
B_r(0)\to B_r(0),\quad x\mto g^n(p^nx)
\]
converge uniformly to a function $f\!: B_r(0)\to B_r(0)$.
Then $f$ admits a $k$-th order expansion as
described in Lemma~{\rm \ref{alongcurv}\/}, and
thus $f\circ \eta$ is $C^k_\K$
for each $C^k_\K$-curve $\eta\!: \K\supseteq W\to B_r(0)$.
\end{prop}
\begin{cor}\label{crexp}
If $G$ is $SC^{k_+}_\K$,
then $\exp$ from Proposition~{\rm \ref{powers}\,(r)\/}
admits a $j$-th order expansion
$($over the ground field $\K)$,
for each integer $j\leq k$,
and $\,\exp\circ \eta$ is $C^k_\K$ for each
$C^k_\K$-curve $\eta\!: \K\supseteq W\to V_r$.
In particular, $\eta_x\!: \Z_p\to V_r$, $z\mto x^z$
and $\Z_p\to V_r$, $z\mto \exp(zx)$
are $C^k_{\Q_p}$, for each $x\in V_r$,
and so is any continuous homomorphism
$\eta\!: \Z_p\to G$.
\end{cor}
\begin{proof}
By the Inverse Function Theorem
for $SC^k$-maps
\cite[Thm.\,7.3]{IMP}, the map
$g:=(\tau_p)^{-1}\!: B_{p^{-1}r}(0)\to B_r(0)$
is $SC^{k_+}_\K$ and hence of class $C^{k_+}_\K$.
The first and second assertion are therefore obvious
from (\ref{explimit}) and Proposition~\ref{difficult}.
For each $x\in V_r$, the map
$\A\to V_r$, $z\mto \exp(zx)$ is $C^k_\K$
(and thus $C^k_{\Q_p}$), being a composition of $\exp$ and
the $C^k_\K$-curve $z\mto zx$.
As $\eta_x(z)=\exp(z\log(x))$ for $z\in \Z_p$,
also $\eta_x$ is $C^k_{\Q_p}$.

Finally, if $\eta\!: \Z_p\to G$ is a continuous homomorphism,
then there exists $N\in \N$ such that
$\eta(z)\in U$ for all $z\in \Z_p$ such that
$|z|_p\leq p^{-N}$. Define $x:=\phi(\eta(p^N))\in V_r$.
Then $z\mto \phi(\eta(p^N z))$
and $\eta_x$ are continuous homomorphisms
$\Z_p\to V_r$ which agree at $1$,
thus on $\langle 1\rangle=\Z$, and hence on all of $\Z_p$,
by continuity. As a consequence, $\eta$ is $C^k_{\Q_p}$ on $p^N\Z_p$
and hence $C^k_{\Q_p}$ on all of $\Z_p$,
by Lemma~\ref{higherdhom}.\vspace{-1mm}
\end{proof}
If $\K$ is locally compact and $G$
finite-dimensional, then the $C^k_\K$-Lie group~$G$ is
also $SC^k_\K$. Otherwise
(when $k_+=k+1$),
we might assume that $G$ is
$C^{k+2}_\K$ to ensure that $G$ is $SC^{k+1}_\K$.\\[2.5mm]
{\bf Proof of Proposition~\ref{difficult}.}
In view of the uniqueness assertion
in Lemma~\ref{alongcurv}, it
suffices to show that every $x_0\in B_r(0)$
has an open neighbourhood $U_{x_0}$
on which $f$
admits a $k$-th order expansion
(by uniqueness, the individual $a_j$'s on the sets $U_{x_0}\times E$
then combine to a well-defined map $a_j$ on $B_r(0)\times E$).\\[3mm]
Thus, fix $x_0\in B_r(0)$ now.
The proof proceeds in two stages.
First, we
show that certain limits exist,
and define certain continuous
maps $b_j$ and $a_j$ using these limits
(Lemma~\ref{yetanother}\,(a), Eqn.\,(\ref{btoa})).
The maps $a_j$
are the natural candidates
for the coefficients
of a $k$-th order expansion for $f$,
and $b_j(x,\sbull)$ is the symmetric $j$-linear
map corresponding to the homogeneous
polynomial $a_j(x,\sbull)$.
The second step, then, will be to show
that these coefficients $a_j$
can be used to expand $f$,
with a remainder~$R$
vanishing of order~$k$
(Lemma~\ref{mischief}).
The decisive tool for the proof
is a notational formalism
(set up in the proof of Lemma~\ref{mischief})
allowing us to track each individual
term in the sum obtained by substituting
$n$ times the $k$-th order Taylor expansion
of $g$ into itself to express $g^n$
(both the multilinear terms
and each individual contribution
to the remainder
term).\footnote{Of course, based
on the symmetry of higher differentials,
many of these terms coincide and could be combined
in one term,
but this does not seem useful
in the present context and would only
make the presentation longer and
more complicated.}
We now begin the proof
with a counterpart
to
the labelling of the contributions
to the $k$-th order Taylor
expansion of $g^n$ just mentioned:
We describe a corresponding notational scheme
allowing us to label
certain compositions of
the differentials of $g$,
which are then used to define the $b_j$'s and $a_j$'s.\\[3mm]
Given $n\in \N$,
consider an $n$-tuple
$(s_1,\ldots, s_n)$,
where, for each $\nu \in\{1,\ldots, n\}$,
\[
s_\nu \!: \{1,\ldots, m_\nu\}\to \{1,2,\ldots,k\}
\]
is a mapping on $\{1,\ldots, m_\nu\}$
for some $m_\nu \in \N$,
such that
\[
m_1=1\qquad\mbox{and}\qquad
m_\nu=\sum_{i=1}^{m_{\nu-1}}s_{\nu-1}(i)\quad\mbox{for $\nu \in\{2,\ldots,n\}$.}
\]
Abbreviate $m_{n+1}:=\sum_{i=1}^{m_n}s_n(i)$.
Given $x\in B_r(0)$, we recursively
define continuous $m_{\nu+1}$-linear maps\vspace{-3mm}
\[
g_{(s_1,\ldots,s_\nu)}^x\!: E^{m_{\nu+1}}\to E
\]
for $\nu =1,\ldots, n$
using higher differentials of $g$ via
\begin{eqnarray*}
g^x_{(s_1)}&:=& \frac{1}{s_1(1)!}d^{s_1(1)}g(f(px),\sbull)\,;\\
g^x_{(s_1,\ldots,s_\nu)}&:=& g^x_{(s_1,\ldots,s_{\nu-1})}\circ
\prod_{i=1}^{m_\nu} \frac{1}{s_\nu(i)!}d^{s_\nu(i)}g(f(p^\nu x),\sbull)\,.
\end{eqnarray*}
For example, if $s_1(1):=2$, $s_2(1):=1$,
and $s_2(2):=3$, then $m_2=2$, $m_3=4$,
\begin{eqnarray*}
g_{(s_1)}^x(u_1,u_2)&=&\frac{1}{2}d^2g\big(f(px),u_1,u_2\big),
\qquad\mbox{and}\\
g_{(s_1,s_2)}^x(u_1,u_2,u_3,u_4)&=&\frac{1}{2}d^2g\Bigl(f(px),\,
dg(f(p^2x),u_1),\, \frac{1}{3!}
d^3g(f(p^2x),u_2,u_3,u_4)\Bigr)
\end{eqnarray*}
for all $x\in B_r(0)$ and $u_1,u_2,u_3,u_4\in E$.\\[3mm]
We let $S_n$ be the set of all $s=(s_1,\ldots,s_n)$ as before
(with variable $m_1,\ldots,m_{n+1}$ depending on~$s$),
and $S_{n,j}:=\{(s_1,\ldots,s_n)\in S_n\!: \sum_{i=1}^{m_n}s_n(i)=j\}$,
for each $j \in \{1,\ldots, k\}$.
For $n\in \N$ and $j\in \{1,\ldots, k\}$, we define
\[
h_{j,n}\!: B_r(0)\times E^j\to E,\quad
h_{j,n}(x,u_1,\ldots, u_j):=\sum_{s\in S_{n,j}}g^x_s
(p^nu_1,\ldots, p^nu_j)\,.
\]
Then $h_{j,n}$ is continuous, and $h_{j,n}(x,\sbull)\!: E^j\to E$
is $j$-linear.\\[3mm]
The mappings $B_r(0)\to \cL^j(E,E)$, $x\mto d^jg(x,\sbull)$
being continuous for each $j\in \{1,\ldots, k\}$
(Lemma~\ref{smallo}),
we find
$\rho_0\in \;]0,r]$
such that
$C_j:=\sup\{\|(j!)^{-1}d^jg(x,\sbull)\|\!: x\in B_{\rho_0}(0)\}<\infty$
for each $j\in \{1,\ldots, k\}$.
There is $N_0\in \N$ such that $p^{-N_0}r\leq \rho_0$
and thus $f(p^nx_0)\in B_{\rho_0}(0)$
for all $n> N_0$.
For each
$n\in \{1,\ldots,N_0\}$, we find $\rho_n\in \;]0,r]$
such that $C_{n,j}:=\sup\{\|(j!)^{-1}d^jg(x,\sbull)\|\!: x\in B_{\rho_n}(f(p^nx_0)\}<\infty$
for each $j\in \{1,\ldots, k\}$.
There is $\rho\in \;]0,r]$ such that
$f(p^nx)\in B_{\rho_n}(f(p^nx_0))$
for all $x\in B_\rho(x_0)$ and $n\in \{1,\ldots, N_0\}$.
Choosing $\rho\leq \rho_0$,
we can achieve that furthermore
$f(p^nx)\in B_{\rho_0}(0)$ for all $n>N_0$
and all $x\in B_\rho(x_0)$, using that $f$ is isometric and $B_{\rho_0}(0)
=B_{\rho_0}(f(p^nx_0))$.\\[3mm]
After replacing
the given norm $\|.\|$
with a suitable positive multiple (and adapting $r$,
$\rho$ and the $\rho_n$'s
accordingly), we may assume that
$C_j\leq 1$ and $C_{n,j}\leq 1$, for all $j\geq 2$
and $n\in \{1,\ldots, N_0\}$.
Then
\begin{equation}\label{actall}
\|(j!)^{-1}
d^jg(f(p^nx),\sbull)\|\leq 1,\quad \mbox{for all $\,j\in \{2,\ldots, k\}$,
$n\in \N$ and $\,x\in B_\rho(x_0)$.}
\end{equation}
It is important to estimate the norms
of the multilinear maps $g^x_s$.
Given $s\in S_{n,j}$, where $j\geq 2$,
there exists
a largest integer $\ell_s\in \{1,\ldots, n\}$
such that $s_{\ell_s}(i)>1$ for some $i\in \{1,\ldots,m_{\ell_s}\}$
(with notation as above). We now show by induction
on $n\in \N$:
\begin{la}\label{strangeind}
$\|g^x_s(p^nu_1,\ldots,p^nu_j)\|\leq p^{-\ell_s}\|u_1\|\cdots\|u_j\|$
holds, for all
$j\in \{2,\ldots,k\}$,
$x\in B_\rho(x_0)$, $u_1,\ldots, u_j\in E$,
and $s\in S_{n,j}$.
\end{la}
\begin{proof}
If $n=1$, then $s_1(1)=j$, $\ell_s=1$. Using (\ref{actall}), we get
\[
\|g^x_s(pu_1,\ldots, pu_j)\|\!=\!
\| (j!)^{-1} d^jg(f(px),p u_1, \ldots, pu_j)\|
\!\leq\! p^{-j} \|u_1\|\!\cdots\!\|u_j\|
\!\leq\! p^{-\ell_s} \|u_1\|\!\cdots\!\|u_j\|.
\]
Now let $n\geq 2$ and suppose the
lemma has been proven up to $n-1$. We have
\[
g^x_s(p^nu_1,\ldots,p^nu_j)
=g^x_{(s_1,\ldots, s_{\ell_s})}(p^{\ell_s}v_1,\ldots,
p^{\ell_s} v_j)
\]
with $v_i:=pg'(f(p^{\ell_s+1}))\cdots pg'(f(p^nx)).u_i$
of norm $\|v_i\|\leq \|u_i\|$, for $i=1,\ldots, j$.

{\em Special case\/}: If $\ell_s=1$ or $\ell_s\geq 2$ and
$(s_1,\ldots,s_{\ell_s-1})\in S_{\ell_s-1,1}$,
then indeed
\begin{eqnarray*}
\|g^x_s(p^n u_1,\ldots,p^n u_j)\|\!\!&\!=\!&\!\! \|
g'(f(px))\cdots g'(f(p^{\ell_s-1}x))(j!)^{-1}
d^jg(f(p^{\ell_s}x),p^{\ell_s}v_1,
\ldots,p^{\ell_s}v_j)\|\\
\!\!&\!\leq\! &\!\! p^{\ell_s-1}p^{-j\ell_s}\|v_1\|\cdots\|v_j\|
\!\leq\!  p^{-(j-1)\ell_s}\|u_1\|\cdots\|u_j\|
\!\leq\!  p^{-\ell_s}\|u_1\|\cdots\|u_j\|.
\end{eqnarray*}
If we are {\em not\/} in the special situation
just described, then $\ell_s\geq 2$
and $t:=(s_1,\ldots, s_{\ell_s-1})\in S_{\ell_s-1,j'}$
with $j':=m_{\ell_s}\in \{2,\ldots, j-1\}$. Hence
\begin{eqnarray*}
\|g^x_s(p^nu_1,\ldots,p^n u_j)\| &\leq&
\|g^x_t\|\cdot \left(\prod_{i=1}^{j'}
\|(s_{\ell_s}(i)!)^{-1}d^{s_{\ell_s}(i)}g(f(p^{\ell_s}x),\sbull)
\|\right)
\|p^{\ell_s}v_1\|\ldots\|p^{\ell_s}v_j\|\\
&\leq & \|g^x_t\| p^{j'}p^{-j\ell_s}\|u_1\|\cdots\|u_j\|
\leq p^{-\ell_t}p^{j'(\ell_s-1)} p^{j'}p^{-j\ell_s}\|u_1\|\cdots\|u_j\|\\
&\leq & p^{-\ell_t-(j-j')\ell_s}\|u_1\|\cdots
\|u_j\|\leq p^{-\ell_s}\|u_1\|\cdots\|u_j\|\,.
\end{eqnarray*}
Here, passing to the second line
we used (\ref{leq1}) and (\ref{actall}).
For the next
inequality, we used that $\|g^x_t\|\leq p^{-\ell_t}p^{j'(\ell_s-1)}$,
by the induction hypothesis. The lemma is established.
\end{proof}
From Lemma~\ref{strangeind}
(case $j\geq 2$)
and (\ref{leq1}) (used if $j=1$), we deduce:
\begin{equation}\label{normcontrol}
\|g^x_s(p^n\sbull,\ldots,p^n\sbull)\|\leq 1,\quad
\mbox{for all $x\in B_\rho(x_0)$, $n\in \N$,
$j\in \{1,\ldots,k\}$, and $s\in S_{n,j}$.}
\end{equation}
Hence also
\begin{equation}\label{dagg}
\|h_{j,n}(x,\sbull)\|\leq 1\quad \mbox{for all
$x\in B_\rho(x_0)$, $j\in \{1,\ldots, k\}$
and $n\in \N$,}
\end{equation}
in view of the ultrametric inequality.
\begin{la}\label{yetanother}
For each $j\in \{1,\ldots, k\}$, we have:
\begin{itemize}
\item[\rm (a)]
For every $x\in B_\rho(x_0)$ and $u_1,\ldots, u_j\in E$, the limit
\[
b_j(x,u_1,\ldots,u_j):=\lim_{n\to\infty} h_{j,n}(x,u_1,\ldots,u_j)
\]
exists in~$E$.
\item[\rm (b)]
For each $x\in B_\rho(x_0)$,
the map $b_j(x,\sbull)\!: E^j\to E$ is
$j$-linear, of norm $\|b_j(x,\sbull)\|\leq 1$.
\item[\rm (c)]
For every $\ve>0$, there exists
%$\delta\in \;]0,\rho]$
$N \in \N$ such that
\[
\|b_j(x,\sbull)-h_{j,n}(x,\sbull)\|\leq \ve\qquad
\mbox{for all $x\in B_\rho(x_0)$ and $n\geq N$.}
\]
\item[\rm (d)]
The map $b_j\!: B_\rho(x_0)\times E^j\to E$
is continuous.
\end{itemize}
\end{la}
\begin{proof}
We observe first that (d) holds if so does (c).
Indeed, given $R>0$ and $\ve'>0$,
apply (c) with $\ve:=\ve'/R^j$.
Then, for every $x\in B_\rho(x_0)$,
$u_1,\ldots, u_j\in
B_R(0)$ and $n\geq N$,
we have
$\|b_j(x,u_1,\ldots,u_j)-h_{j,n}(x,u_1,\ldots, u_j)\|\leq \ve R^j=\ve'$,
showing that $h_{j,n}\to b_j$
uniformly on $B_\rho(x_0)\times B_R(0)^j$.
Hence $b_j$ is continuous on
$B_\rho(x_0)\times B_R(0)^j$
for each $R>0$ and hence continuous on all of $B_\rho(x_0)\times E^j$.
It remains to prove (a)--(c).\vspace{1mm}

{\em The case $j=1$.}
For each $n\in \N$, the set
$S_{n,1}$ has one element $s$
only, and we have
\begin{equation}\label{leq12}
h_{1,n}(x,\sbull)=g^x_s(p^n\sbull)=pg'(f(px))\circ\cdots \circ p g'(f(p^nx))\quad
\mbox{for each $x\in B_\rho(x_0)$.}
\end{equation}
Given $\ve>0$, say $\ve\in \;]0,1[$ without loss
of generality, there exists $\sigma\in \;]0,r]$
such that
\begin{equation}\label{sigm}
\|pg'(y)-\id_E\|<\ve \quad\mbox{for all $y\in B_\sigma(0)\sub E$.}
\end{equation}
Choose $N\in \N$ such that $p^{-N}r\leq \sigma$;
then $f(p^nx)\in B_\sigma(0)$ for all $x\in B_r(0)$
and all $n\geq N$.
For all $x\in B_\rho(x_0)$ and $n_0,n_1\geq N$,
where $n_1\geq n_0$ without loss of generality,
we have
\[
h_{1,n_0}(x,\sbull)^{-1}\circ
h_{1,n_1}(x,\sbull)=pg'(f(p^{n_0+1}x))\circ \cdots
\circ pg'(f(p^{n_1}x))=\id_E +A
\]
for some $A\in \cL(E)$ with $\|A\|<\ve$,
using that $B_\ve^{\cL(E)}(\id_E)$ is a subgroup
of $\cL(E)^\times$.
Hence
\begin{eqnarray}
\|h_{1,n_1}(x,\sbull)-h_{1,n_0}(x,\sbull)\| & = &\|h_{1,n_0}(x,\sbull)\circ A\|
\leq \|h_{1,n_0}(x,\sbull)\|\cdot \|A\|\nonumber\\
&\leq& \|A\|<\ve
\quad\mbox{for all $x\in B_\rho(x_0)$ and $n_0,n_1\geq N$,}\label{unifo}
\end{eqnarray}
using (\ref{dagg}).
Thus $(h_{1,n}(x,\sbull))_{n\in \N}$ is a Cauchy sequence
in $\cL(E)$ for $x\in B_\rho(x_0)$,
and hence converges to some $b_1(x,\sbull)\in \cL(E)$.
Since $\|h_{1,n}(x,\sbull)\|\leq 1$ by (\ref{dagg}),
we also have $\|b_1(x,\sbull)\|\leq 1$,
whence (a) and (b) hold for $j=1$. Taking $n:=n_0$
and letting $n_1\to \infty$ in
(\ref{unifo}), we see that
$\|b_1(x,\sbull)-h_{1,n}(x,\sbull)\|\leq\ve$
for all $x\in B_\rho(x_0)$ and $n\geq N$,
whence (c) holds.\vspace{1mm}

{\em Now let $j\geq 2$.}
If we can show that,
for every $\ve'>0$, there is $N\in \N$
such that
\begin{equation}\label{firstform}
\|h_{j,n_1}(x,\sbull)-h_{j,n_0}(x,\sbull)\|\leq\ve'
\quad\mbox{for all $x\in B_\rho(x_0)$ and $n_0,n_1\geq N$,}
\end{equation}
then
(a), (b) and (c) follow exactly as in the case $j=1$.
To establish (\ref{firstform}), and hence to
complete the proof of the lemma,
it suffices to find $N\in \N_0$ such that
\begin{equation}\label{secform}
\left\{
\begin{array}{l}
h_{j,n_1}(x,u_1,\ldots, u_j)-h_{j,n_0}(x,u_1,\ldots,u_j)\in \wb{B}_\ve(0)\\
\mbox{for all $x\in B_\rho(x_0)$, $n_0,n_1\geq N$, and $u_1,\ldots,u_j\in E$
of norm $\leq 1$,}
\end{array}
\right.
\end{equation}
where $\ve:=p^{-j}\ve'$
(cf.\ \cite[p.\,59]{Roo}).
To this end, assuming w.l.o.g.\
$\ve<1$,
we choose $\sigma>0$
and $N\in \N$ as in the proof
of the case $j=1$, with $p^{-N}<\ve$ now.
To see that (\ref{secform})
holds with this~$N$,
let $x\in B_\rho(x_0)$,
and $u_1,\ldots,u_j\in E$ be vectors of norm $\leq 1$.
For $n\geq N$ and $s\in S_{n,j}$
with $\ell_s>N$,
we have
\[
\|g^x_s(p^nu_1,\ldots, p^nu_j)\|\leq p^{-N}\|u_1\|\cdots\|u_j\|\leq\ve\,,
\]
by Lemma~\ref{strangeind}.
Hence, using that $\wb{B}_\ve(0)$ is an additive subgroup of~$E$,
we obtain
\begin{equation}\label{step1}
h_{j,n}(x,u_1,\ldots,u_n)\in \sum_{s\in S_{n,j,N}}g^x_s(p^nu_1,\ldots,
p^nu_j)\;+ \wb{B}_\ve(0)\,,
\end{equation}
where $S_{n,j,N}:=\{s\in S_{n,j}\!: \ell_s\leq N\}$.
Abbreviating
\[
S_{\ell,j}^*:=\{s\in S_{\ell,j}\!: \ell_s=\ell\}
\quad\mbox{and}\quad A_{\ell,n}:=pg'(f(p^{\ell+1}x))\circ
\cdots\circ pg'(f(p^nx))\in \cL(E)
\]
for $\ell\in \{1,\ldots, N\}$
(where $A_{\ell,n}=\id_E$ if $\ell=n$),
the above sum can be re-written as
\begin{equation}\label{step2}
\sum_{s\in S_{n,j,N}}g^x_s(p^nu_1,\ldots,
p^nu_j)
=\sum_{\ell=1}^N\sum_{s\in S_{\ell,j}^*}
g_s^x(p^\ell A_{\ell,n}u_1,\ldots , p^\ell A_{\ell,n}u_j)\,.
\end{equation}
If $n_1\geq n_0\geq N$, $\ell\in \{1,\ldots,N\}$, and $s\in S_{\ell,j}^*$,
then
$B:=pg'(f(p^{n_0+1}x))\circ \cdots\circ pg'(f(p^{n_1}x))
\in B^{\cL(E)}_\ve(\id_E)$  (cf.\ proof of (\ref{unifo}))
and thus
\begin{eqnarray}
\lefteqn{g_s^x(p^\ell A_{\ell,n_1}u_1,\ldots p^\ell A_{\ell,n_1}u_j)-
g_s^x(p^\ell A_{\ell,n_0}u_1,\ldots p^\ell A_{\ell,n_0}u_j)}\nonumber\\
&=&
g_s^x(p^\ell A_{\ell,n_0}Bu_1,\ldots p^\ell A_{\ell,n_0}Bu_j)-
g_s^x(p^\ell A_{\ell,n_0}u_1,\ldots p^\ell A_{\ell,n_0}u_j)
\in \wb{B}_\ve(0)\label{step3}
\end{eqnarray}
by Lemma~\ref{multilin}, in view of (\ref{normcontrol}),
$\|A_{\ell,n_0}\|\leq 1$,
$\|B\|\leq 1$,
and $\|B-\id_E\|<\ve$.
Combining (\ref{step1}),
(\ref{step2}) and (\ref{step3}),
we see that
\begin{eqnarray*}
\lefteqn{h_{j,n_1}(x,u_1,\ldots,u_j)-h_{j,n_0}(x,u_1,\ldots,u_j)}\\
&\in &
\sum_{\ell=1}^N\sum_{s\in S_{\ell,j}^*}
\left(g_s^x(p^\ell A_{\ell,n_1}u_1,\ldots p^\ell A_{\ell,n_1}u_j)
- g_s^x(p^\ell A_{\ell,n_0}u_1,\ldots p^\ell A_{\ell,n_0}u_j)
\right)+\wb{B}_\ve(0) =\wb{B}_\ve(0)\,,
\end{eqnarray*}
whence (\ref{secform}) holds.
This completes the proof.
\end{proof}
For each $j\in \{1,\ldots, k\}$, we define
$a_j\!: B_\rho(x_0)\times E\to E$ via
\begin{equation}\label{btoa}
a_j(x,u):=b_j(x,u,\ldots,u)\qquad\mbox{for $x\in B_\rho(x_0)$, $u\in E$.}
\end{equation}
Then $a_j$ is continuous, and $a_j(x,\sbull)$
is a homogeneous polynomial of degree~$j$,
for each $x\in B_\rho(x_0)$.
To complete the proof of Proposition~\ref{difficult},
we show:
\begin{la}\label{mischief}
$R\!: B_\rho(x_0)\times B_\rho(x_0)\to E$, $R(x,y):=
f(y)-f(x)-\sum_{j=1}^ka_j(x,y-x)\,$
is a continuous mapping such that
\[
f(y)=f(x)+\sum_{j=1}^k a_j(x,y-x)+R(x,y)\quad\mbox{for all
$\,x,y\in B_\rho(x_0)$,}
\]
$R(x,x)=0$ for all $\,x\in B_\rho(x_0)$,
and
\begin{equation}\label{vani}
\lim_{(x,y)\to (z_0,z_0)}\frac{\|R(x,y)\|}{\|x-y\|^k}=0
\qquad\mbox{for all $\,z_0\in B_\rho(x_0)\,$ $($where $x\not=y)$.}
\end{equation}
\end{la}
\begin{proof}
First, we set up a notational formalism
enabling us to explicitly label each term in the sum
obtained by expanding every factor $g$ of $g^n$
into its $k$-th order Taylor expansion:
\begin{equation}\label{formalism}
g^n(p^ny)-g^n(p^nx)=
\sum_{s\in S^0_n}H^{x,y}_{n,s}(p^n(y-x),\ldots,p^n(y-x))\quad
\mbox{for all $x,y\in B_r(0)$.}
\end{equation}
Let $P\!: B_r(0)\times B_r(0)\to E$,
$P(x,y):=g(y)-g(x)-\sum_{j=1}^k\frac{1}{j!}d^jg(x,y-x,\ldots,y-x)$
be the remainder term of the $k$-th order Taylor expansion
of~$g$. To enable a unified notation
for differentials and remainder terms,
for $x,y\in B_r(0)$ and $j\in \{0,\ldots, k\}$ we define
\[
c_j^{x,y}:=\left\{
\begin{array}{cl}
\frac{1}{j!}d^jg(x,\sbull)\in \cL^j(E,E)
& \mbox{if $\;j\in \{1,\ldots, k\}$;}\\
\!\!\!\!\!\!P(x,y) \in E\;\;\; & \mbox{if $\;j=0$.}
\end{array}
\right.
\]
In the following, expressions like
$c^{x,y}_0(u_i,\ldots,u_j)$ with $i>j$
have to be read as $c^{x,y}_0$; they denote elements of~$E$
(not functions).
Given $n\in \N$, we let $S^0_n$ be the set of all $(s_1,\ldots,s_n)$
where $s_\nu \!: \{1,\ldots, m_\nu\}\to \{0,1,\ldots,k\}$
for $\nu =1,\ldots, n$
for certain $m_\nu\in \N_0$ such that
$m_1=1$ and
$\sum_{i=1}^{m_\nu}s_\nu(i)=m_{\nu+1}$ for $\nu=1,\ldots, n-1$;
set $m_{n+1}:=\sum_{i=1}^{m_n}s_n(i)$.
Given $x,y\in B_r(0)$, $n\in \N$, and $s=(s_1,\ldots,s_\nu)\in S_\nu^0$, where
$\nu\in \N$ with $\nu \leq n$,
we define
\[
H^{x,y}_{n,s}(u_1,\ldots,u_{m_{\nu+1}})
:=c^{g^{n-1}(p^nx),g^{n-1}(p^ny)}_{s_1(1)}(u_1,\ldots, u_{m_{\nu+1}})
\quad\mbox{for $u_1,\ldots,u_{m_{\nu+1}}\in E$}
\]
if $\nu=1$, and recursively for $\nu\geq 2$
\[
\begin{array}{l}
H^{x,y}_{n,s}(u_1,\ldots,u_{m_{\nu+1}})\;:=\\
\; H^{x,y}_{n,(s_1,\ldots,s_{\nu-1})}
\bigl(
c^{g^{n-\nu}(p^nx),g^{n-\nu}(p^ny)}_{s_\nu(1)}(u_1,\ldots,u_{s_\nu(1)}),\,
c^{g^{n-\nu}(p^nx),g^{n-\nu}(p^ny)}_{s_\nu(2)}(u_{s_\nu(1)+1},
\ldots,u_{s_\nu(1)+s_\nu(2)}),\\
\qquad\qquad \qquad\ldots,
c^{g^{n-\nu}(p^nx),g^{n-\nu}(p^ny)}_{s_\nu(m_\nu)}(u_{m_{\nu+1}-s_\nu(m_\nu)+1},\ldots,u_{m_{\nu+1}})\bigr)\,.
\end{array}
\]
For example, if $n=2$, $m_1=1$, $m_2=3$ and $s:=(s_1,s_2)$ with
$s_1(1):=3$, $s_2(1):=1$, $s_2(2)=0$ and $s_2(3):=2$,
then $m_3=3$ and
\[
H^{x,y}_{n,s}(u_1,u_2,u_3)=\frac{1}{3!}\, d^3g\Bigl( g(p^2x),\,
dg(p^2x,u_1),\, P(p^2x,p^2y),\,\frac{1}{2}d^2g(p^2x, u_2,u_3)\Bigr)\,.
\]
Using the notational formalism,
when $k=n=2$ we calculate
\begin{eqnarray*}
\lefteqn{g^2(p^2y)-g^2(p^2x)=g(g(p^2y))-g(g(p^2x))}\\
&=&
g'(g(p^2x)).(g(p^2y)-g(p^2x))+\frac{1}{2}
d^2g(g(p^2x),g(p^2y)-g(p^2x),
g(p^2y)-g(p^2x))\\
& & \qquad + \, P(g(p^2x),g(p^2y))\\
&=&
\sum_{s\in S_2^0}H_{2,s}^{x,y}(p^2(y-x),\ldots, p^2(y-x))\,;\vspace{-4mm}
\end{eqnarray*}
here, we used the Taylor expansion of~$g$ around $g(p^2x)$ to pass to the
second line, and then re-wrote
$g(p^2y)-g(p^2x)$ using the Taylor expansion of $g$ around $p^2x$,
to pass to the third
(the reader might
write down all 13 summands to check this).
Likewise, expanding each of the $n$ factors of $g^n$
in turn,
a moment's reflection shows that (\ref{formalism})
holds, for all $n\in \N$.\\[3mm]
We now fix $\ve\in \;]0,1[$ and $z_0\in B_\rho(x_0)$
for the rest of the proof.
Keeping the notation set up before
Lemma~\ref{strangeind},
we find $\delta_0\in \;]0,\rho_0]$
such that $\|P(x,y)\|\leq \ve\|x-y\|^k$
for all $x,y\in B_{\delta_0}(0)$.
There is $N_1\geq N_0$ such that $p^{-N_1}r\leq \delta_0$.
For each $n\in \{1,\ldots,N_1\}$, we find $\delta_n\in \;]0,r]$
such that $\|P(x,y)\|\leq \ve\|y-x\|^k$
for all $x,y\in B_{\delta_n}(f(p^nz_0))$.
If $n\leq N_0$, we assume that
$\delta_n\leq \rho_n$ here.
Next, we find $\delta\in \;]0,\min\{\rho,\ve\}[$
(whence $\delta<1$ in particular)
such that $f(p^nx)\in B_{\delta_n}(f(p^nz_0))$
for each $n\in \{1,\ldots, N_1\}$
and each $x\in B_\delta(z_0)\sub B_\rho(x_0)$.
Note that $f(p^nx)\in B_{\delta_0}(0)$ for each $n\geq N_1$
and each $x\in B_r(0)$. Our goal is to show that
\begin{equation}\label{finish}
\|R(x,y)\|\leq \ve\,\|y-x\|^k\quad\mbox{for all $x,y\in B_\delta(z_0)$,}
\end{equation}
thus establishing (\ref{vani}).
Since $g^n(p^nx)\to f(x)$ uniformly in~$x$ by hypothesis,
there exists $N_2\geq N_1$ such that
\begin{equation}\label{unf}
\|g^{n-\nu}(p^nx)-f(p^\nu x)\|\leq \delta_\nu\quad
\mbox{for all $x\in B_r(0)$, $\nu\in \{1,\ldots,N_1\}$,
and $n\geq N_2\,$,}
\end{equation}
whence $g^{n-\nu}(p^nx)\in B_{\delta_\nu}(f(p^\nu z_0))$
if $x\in B_\delta(z_0)$.
On the other hand,
for every $n\geq N_2$ and $\nu\in \N$ such that $N_1<\nu\leq n$,
we have $\|g^{n-\nu}(p^nx)\|=p^{-\nu}\|x\|<p^{-\nu}r<p^{-N_1}r
\leq \delta_0$.
As a consequence,
we have
\begin{equation}\label{summry}
\|P(g^{n-\nu}(p^nx),g^{n-\nu}(p^ny))\| \leq \ve
\|g^{n-\nu}(p^nx)-g^{n-\nu}(p^ny)\|^k
\leq \ve p^{-\nu k} \|y-x\|^k
\end{equation}
for all
$x,y\in B_\delta(z_0)$, $n\geq N_2$ and $\nu\in \{1,\ldots, n\}$.
Let
$x\in B_\delta(z_0)$, $n\geq N_2$, and $\nu\in \{1,\ldots,n\}$.
If $\nu\leq N_0$,
then $g^{n-\nu}(p^nx)\in B_{\delta_\nu}(f(p^\nu z_0))
\sub B_{\rho_\nu}(f(p^\nu z_0))\sub B_{\rho_\nu}(f(p^\nu x_0))\,$;
if $\nu>N_0$, then $\|g^{n-\nu}(p^nx)\|=p^{-\nu}\|x\|<p^{-N_0}r\leq \rho_0$
and thus $g^{n-\nu}(p^nx)\in B_{\rho_0}(0)$. Hence
\begin{equation}\label{dag2}
\|(j!)^{-1}d^jg(g^{n-\nu}(p^nx),\sbull)\|\leq 1\quad
\mbox{for all $x\in B_\delta(z_0)$,
$2\leq j\leq k$, $n\geq N_2$,
and $1\leq \nu \leq n$.}
\end{equation}
To establish (\ref{finish}), we prove estimates
on the norms $\|H^{x,y}_{n,s}\|$,
which will enable us to get rid of all summands
involving remainder terms,
or which are multilinear of order exceeding $k$:\\[3mm]
{\bf Claim~1.} {\em For any $n\geq N_2$, $\nu\in \{1,\ldots, n\}$,
$x,y\in B_\delta(z_0)$, and $s\in S^0_\nu$, we
have
\begin{equation}\label{weaker}
\|H^{x,y}_{n,s}\|\leq p^{\nu \, m_{\nu+1}}\,.
\end{equation}
If $s\in S_\nu^0\setminus S_\nu$ here, then
furthermore}
\begin{equation}\label{stronger}
\|H^{x,y}_{n,s}\|\leq \ve\|y-x\|^k p^{\nu \, m_{\nu+1}}\,.
\end{equation}
{\bf Claim~2.} {\em For any $n\geq N_2$,
$x,y\in B_\delta(z_0)$ and $s\in S^0_n$
with $s\not\in S_n$ or $m_{n+1}>k$,
we have}
\begin{equation}\label{gtrid}
\|H^{x,y}_{n,s}(p^n(x-y),\ldots,p^n(y-x))\|\leq \ve \|y-x\|^k\,.
\end{equation}
Once these claims are proved, for $n\geq N_2$ and $x,y\in B_\delta(z_0)$
we shall simply have
\begin{eqnarray}
g^n(p^ny)-g^n(p^nx)&=&\sum_{s\in S_n^0}H^{x,y}_{n,s}(p^n(y-x),\ldots,p^n(y-x))\nonumber\\
&\in&
\sum_{j=1}^k \wt{h}_{j,n}(x,y-x,\ldots,y-x)\,
+\, \wb{B}_{\ve\|y-x\|^k}(0)\label{simplify}
\end{eqnarray}
with $\wt{g}^x_s:=H^{x,x}_{n,s}$
for $s\in S_n$ and $\wt{h}_{j,n}(x,u_1,\ldots, u_j):=\sum_{s\in S_{n,j}}
\wt{g}^x_s(p^nu_1,\ldots, p^nu_j)$.\\[3mm]
{\bf Proof of Claim~1.}
Fix $n\geq N_2\,$; the proof is by induction on $\nu\in \{1,\ldots, n\}$.
Assume $\nu=1$ first;
thus $s=(s_1)$. If $s_1(1)=0$,
then $\|H_{n,s}^{x,y}\|=\|P(g^{n-1}(p^nx),g^{n-1}(p^ny))\|\leq
\ve p^{-k}\|y-x\|^k\leq \ve \|y-x\|^k$ by (\ref{summry}),
whence (\ref{stronger}) holds and also (\ref{weaker}),
as $\ve\|y-x\|^k\leq 1=p^{m_2}$, because $m_2=0$.
If $j:=s_1(1)\geq 1$,
then
$\|H^{x,y}_{n,s}\|=\|(j!)^{-1}d^jg(g^{n-1}(p^nx),\sbull)\|
\leq p\leq p^{1 m_2}$ by (\ref{leq1}) and (\ref{dag2}),
whence (\ref{weaker}) holds
(we need not check (\ref{stronger}), because $s\in S_1$).

{\em Induction step.} Let $\nu\in \{2,\ldots,n\}$, and
suppose the
assertion is correct for $\nu$ replaced with $\nu-1$.
Given $s\in S^0_\nu$, abbreviate $t:=(s_1,\ldots,s_{\nu-1})$.
There are four cases:\vspace{1.3mm}

{\em Case\/} 1: {\em $m_\nu=0$.} Then also $m_{\nu+1}=0$, and
$H^{x,y}_{n,s}=
H^{x,y}_{n,t}$ with $t\in S^0_{\nu-1}\setminus S_{\nu-1}$
and thus $\|H^{x,y}_{n,s}\|=\|H^{x,y}_{n,t}\|\leq \ve\|y-x\|^kp^0$
by induction, whence (\ref{stronger}) and (\ref{weaker})
hold.\vspace{1.3mm}

{\em Case\/} 2: {\em $m_\nu>0$, $s\in S^0_\nu\setminus
S_\nu$, and $\,t\in S_{\nu-1}$.}
Then $s_\nu^{-1}(\{0\})$ is a non-empty set;
let $j\in \{1,\ldots,m_\nu\}$ be its number of elements.
As $H^{x,y}_{n,t}$
is continuous $m_\nu$-linear of norm
$\|H^{x,y}_{n,t}\|\leq p^{(\nu-1)m_\nu}$
by induction and $H^{x,y}_{n,s}$ is obtained by
inserting $j$ times $P(g^{n-\nu}(p^nx),g^{n-\nu}(p^ny))$
and $(m_\nu\!-\!j)$ times multilinear maps of norms
$\leq p$ into $H^{x,y}_{n,t}$ (see (\ref{leq1}), (\ref{dag2})),
we get
\begin{eqnarray*}
\|H^{x,y}_{n,s}\| &\leq &\|H^{x,y}_{n,t}\|\cdot p^{m_\nu-j}\cdot
\|P(g^{n-\nu}(p^nx),g^{n-\nu}(p^ny))\|^j
\leq p^{(\nu-1)m_\nu}p^{m_\nu-j}
\ve^jp^{-\nu k j}\|y-x\|^{kj}\\
&\leq & p^{\nu (m_\nu-j)}
\ve \|y-x\|^k
\leq p^{\nu \, m_{\nu+1}}
\ve \|y-x\|^k\,,
\end{eqnarray*}
using that
$\|P(g^{n-\nu}(p^nx)g^{n-\nu}(p^ny))\|\leq \ve p^{-\nu k}\|y-x\|^k$
by (\ref{summry}),
$\ve\leq 1$,
$\|y-x\|\leq 1$,
and $m_{\nu+1}=\sum_{i=1}^{m_\nu}s_\nu(i)
=\sum_{i\in s_\nu^{-1}(\N)} s_\nu(i)
=(m_\nu-j)+\sum_{i\in s_\nu^{-1}(\N)}(s_\nu(i)-1)
\geq m_\nu-j$.
Thus (\ref{stronger}) and
(\ref{weaker}) hold.\vspace{1.3mm}

{\em Case\/} 3: {\em $m_\nu>0$ and $t\in S_{\nu-1}^0\setminus
S_{\nu-1}$}. Then, by induction,
$H^{x,y}_{n,t}$ is continuous $m_\nu$-linear,
of norm $\leq \ve \|y-x\|^kp^{(\nu-1)m_\nu}$.
Let $j\in \{0,\ldots, m_\nu\}$
be the number of zeros of $s_\nu$.
Since $H^{x,y}_{n,s}$
is obtained from $H^{x,y}_{n,t}$
by inserting $(m_\nu\!-\!j)$ times multilinear maps of norms $\leq p$
and $j$ times
the element $P(g^{n-\nu}(p^nx),g^{n-\nu}(p^ny))$
of norm $\leq p^{-\nu}$ (cf.\ (\ref{summry})), we obtain
$\|H^{x,y}_{n,s}\|\leq \|H^{x,y}_{n,t}\|\cdot p^{m_\nu-j}p^{-\nu j}
\leq \ve\|y-x\|^kp^{(\nu-1)m_\nu}p^{m_\nu-j}p^{-\nu j}
\leq \ve \|y-x\|^kp^{\nu (m_\nu-j)}
\leq \ve \|y-x\|^kp^{\nu\, m_{\nu+1}}$.\vspace{1.3mm}

{\em Case\/} 4: {\em $s\in S_\nu$.}
Then $t\in S_{\nu-1}$
and $1\leq m_\nu\leq m_{\nu+1}$.
By induction, $H^{x,y}_{n,t}$ is a continuous $m_\nu$-linear
mapping of norm
$\|H^{x,y}_{n,t}\|\leq p^{(\nu-1) m_\nu}$.
Therefore
$\|H^{x,y}_{n,s}\|\leq
\|H^{x,y}_{n,t}\|\cdot \prod_{i=1}^{m_\nu}
\|(s_\nu(i)!)^{-1}d^{s_\nu(i)}g(g^{n-\nu}(p^nx),\sbull)\|
\leq p^{(\nu-1)m_\nu}p^{m_\nu}=p^{\nu \, m_\nu}
\leq p^{\nu \,m_{\nu+1}}$. This completes the proof
of Claim~1.\vspace{3mm}\Punkt

\noindent
{\bf Proof of Claim~2.}
Let $n\geq N_2$, $x,y\in B_\delta(z_0)$
and $s\in S^0_n$.
If $s\in S^0_n\setminus S_n$,
then indeed $\|H^{x,y}_{n,s}(p^n(x-y),\ldots,p^n(x-y))\|
\leq \ve\|y-x\|^kp^{n m_{n+1}}\|p^n(y-x)\|^{m_{n+1}}
=\ve\|y-x\|^k\|y-x\|^{m_{n+1}}\leq\ve\|y-x\|^k$,
using (\ref{stronger}).
If $s\in S_n$ and $m_{n+1}>k$,
then
$\|H^{x,y}_{n,s}(p^n(y-x),\ldots,p^n(y-x))\|
\leq \ve p^{nm_{n+1}}\|p^n(y-x)\|^{m_{n+1}}
=\ve \|y-x\|^{m_{n+1}-k}\|y-x\|^k
\leq \ve\|y-x\|^k$, using (\ref{weaker}).\vspace{3mm}\Punkt

\noindent
To prove (\ref{finish}), we now
fix $x,y\in B_\delta(z_0)$ for the rest of the proof,
assuming without loss of generality that $x\not=y$
(the omitted case being trivial).
Thus $\wb{\ve}:=\ve\,\|y-x\|^k>0$.
There exists $N_3\geq N_2$
such that $\|f(x)-g^n(p^nx)\|\leq\wb{\ve}$,
$\|f(y)-g^n(p^ny)\|\leq\wb{\ve}$,
and
\begin{equation}\label{sparkle}
\|b_j(x,y-x,\ldots,y-x)-h_{j,n}(x,y-x,\ldots,y-x)\|\leq \wb{\ve}\qquad
\mbox{for $j\in \{1,\ldots, k\}$,}
\end{equation}
for all
$n\geq N_3$.
By the preceding and
(\ref{simplify}), we have
\[
f(y)-f(x)\in g^n(p^ny)-g^n(p^nx)+\wb{B}_{\wb{\ve}}(0)
=\sum_{j=1}^k\wt{h}_{j,n}(x,y-x,\ldots,y-x) +\wb{B}_{\wb{\ve}}(0)
\]
for all $n\geq N_3$ and thus
\[
R(x,y)\in \sum_{j=1}^k\bigl(\wt{h}_{j,n}(x,y-x,\ldots,y-x)-h_{j,n}
(x,y-x,\ldots,y-x)\bigr)+ \wb{B}_{\wb{\ve}}(0)\,,
\]
by (\ref{sparkle}).
Hence (\ref{finish}) will hold if we can show that
\begin{equation}\label{creach}
\|h_{j,n}(x,y-x,\ldots,y-x)-\wt{h}_{j,n}(x,y-x,\ldots,y-x)\|\leq \wb{\ve}
\end{equation}
for all $j\in \{1,\ldots, k\}$ and all
sufficiently large~$n$.

{\em Assume $j=1$ first.}
We recall that
$h_{1,n}(x,\sbull)=pg'(f(px))\circ\cdots\circ p g'(f(p^nx))$
and $\wt{h}_{1,n}(x,\sbull)=pg'(g^{n-1}(p^nx))\circ
\cdots\circ pg'(g^0(p^nx))$.
There is $\sigma\in \;]0,r]$ such that $\|pg'(z)-\id_E\|<\wb{\ve}$
for all $z\in B_\sigma(0)$, and $N_4\geq N_3$
such that $p^{-N_4}r<\sigma$.
Then
$\|pg'(f(p^nx))-\id_E\|<\wb{\ve}$
and
$\|pg'(g^{n-\nu}(p^nx))-\id_E\|<\wb{\ve}$
for all $n>N_4$ and $\nu\in \{N_4,\ldots, n\}$.
Recalling (\ref{dagg}), we deduce for each
$n>N_4$ that
\begin{eqnarray*}
h_{1,n}(x,y-x) & = & h_{1,n-1}(x,\, pg'(f(p^nx)).(y-x))\\
&= &h_{1,n-1}(x,y-x)+
h_{1,n-1}(x,\,(pg'(f(p^nx))-\id_E).(y-x)))\\
&\in&
h_{1,n-1}(x,y-x)+\wb{B}_{\wb{\ve}}(0)\,.
\end{eqnarray*}
Repeating this argument, we arrive at
$h_{1,n}(x,y-x)\in h_{1,N_4}(x,y-x)+ \wb{B}_{\wb{\ve}}(0)$
for all $n>N_4$.
Similarly, we see that
$\wt{h}_{1,n}(x,y-x)\in pg'(g^{n-1}(p^nx))\circ\cdots\circ
pg'(g^{n-N_4}(p^nx)).(y-x)+
\wb{B}_{\wb{\ve}}(0)$, for all $n>N_4$.
Hence
\[
h_{1,n}(x,x-y)-\wt{h}_{1,n}(x,y-x)\in
J_n.(y-x)
+
\wb{B}_{\wb{\ve}}(0)\quad
\mbox{for all $n>N_4$, where}
\]
$J_n:=
pg'(f(px))\circ\cdots\circ
pg'(f(p^{N_4}x))\,-\,
pg'(g^{n-1}(p^nx))\circ\cdots\circ
pg'(g^{n-N_4}(p^nx))$.
Because $g^{n-\nu}(p^nx)\to f(p^\nu x)$
as $n\to \infty$, because $g'$ is continuous and
also the composition map $\cL(E)\times \cL(E)\to \cL(E)$
is continuous,
we see that $J_n\to 0$ as $n\to\infty$,
whence there exists $N \geq N_4$
such that $\|J_n\|\leq\wb{\ve}$ for all
$n\geq N$. Thus
(\ref{creach}) holds for $j=1$
and all $n\geq N$.\vspace{1mm}

{\em Now assume that $j\in \{2,\ldots, k\}$.}
Repeating the proof
of (\ref{step1})
with $\wb{\ve}$ instead of $\ve$,
we see that
$h_{j,n}(x,u_1,\ldots,u_j)\in \sum_{s\in S_{n,j,N}}g^x_s(p^nu_1,\ldots,
p^n u_j)+ \wb{B}_{\wb{\ve}}(0)$
for all $j\in \{2,\ldots,k\}$, $n\geq N$,
and $u_1,\ldots,u_j\in E$ of norm $\leq 1$
(in view of our choice of $\sigma$ and $N_4$).
In view of (\ref{dag2}), copying the proof
of Lemma~\ref{strangeind}
we see that $\|\wt{g}^x_s(p^n\sbull,\ldots,p^n\sbull)\|\leq p^{-\ell_s}$
for all $n\geq N_2$, $j\in \{2,\ldots, k\}$ and
$s\in S_{n,j}$,
entailing that
\[
\wt{h}_{j,n}(x,u_1,\ldots,u_j)\in \sum_{s\in S_{n,j,N}}
\wt{g}^x_s(p^nu_1,\ldots,
p^nu_j)+ \wb{B}_{\wb{\ve}}(0)
\]
for all $j\in \{2,\ldots, k\}$,
$n\geq N$, and $u_1,\ldots, u_j\in E$ of norm
$\leq 1$.
Setting
$\wt{A}_{\ell,n}:=pg'(g^{n-\ell-1}(p^nx))\circ \cdots\circ pg'(p^nx)$
for $\ell\in \{1,\ldots, N\}$,
we can write
$\wt{g}^x_s(p^nu_1,\ldots, p^nu_j)
=\wt{g}^x_{n,t}(p^{\ell_s}\wt{A}_{\ell_s,n}u_1,\ldots,
p^{\ell_s}\wt{A}_{\ell_s,n}u_j)$,
where
$t:=(s_1,\ldots, s_{\ell_s})$ and
$\wt{g}^x_{n,t}:=H^{x,x}_{n,t}$.
As a consequence,
$h_{j,n}(x,u_1,\ldots,u_j)-\wt{h}_{j,n}(x,u_1,\ldots,u_j)$ is contained in
\[
\sum_{\ell=1}^N\sum_{s\in S^*_{\ell,j}}\Big(
g_s^x(p^\ell A_{\ell,n}u_1,\ldots , p^\ell A_{\ell,n}u_j)
- \wt{g}_{n,s}^x(p^\ell \wt{A}_{\ell,n}u_1,\ldots , p^\ell \wt{A}_{\ell,n}u_j)\Big)
+\wb{B}_{\wb{\ve}}(0)\,,
\]
with $A_{\ell,n}$ as in the proof of Lemma~\ref{yetanother}.
Henceforth, we fix $u_i:=y-x$ for $i=1,\ldots, k$.
To establish (\ref{creach}),
we now only need to show that there exists $N_5\geq N$ such that
\begin{equation}\label{end-2}
g_s^x(p^\ell A_{\ell,n}u_1,\ldots , p^\ell A_{\ell,n}u_j)
- \wt{g}_{n,s}^x(p^\ell \wt{A}_{\ell,n}u_1,\ldots , p^\ell \wt{A}_{\ell,n}u_j)
\in \wb{B}_{\wb{\ve}}(0)
\end{equation}
for all $j\in \{2,\ldots,k\}$, $\ell\in \{1,\ldots,N\}$,
$s\in S_{\ell,j}^*$,
and $n\geq N_5$.
Recall from the case $j=1$
that $\|pg'(z)-\id_E\|<\wb{\ve}$
for all $z\in E$ such that $|z|\leq p^{-N}r$
(because $N\geq N_4$).
Also recall that
$\|\wt{g}^x_{n,s}(p^\ell\sbull,\ldots,p^\ell\sbull)\|\leq 1$
and $\|pg'(z)\|\leq 1$ for all~$z$.
Replacing each of $pg'(g^0(p^nx)),\ldots, pg'(g^{n-N-1}(p^nx))$
with $\id_E$ in turn,
we therefore obtain with Lemma~\ref{multilin} (applied $n-N$ times)
that
\begin{equation}\label{end-1}
\wt{g}_{n,s}^x(p^\ell \wt{A}_{\ell,n}u_1,\ldots , p^\ell \wt{A}_{\ell,n}u_j)
\in
\wt{g}_{n,s}^x(p^\ell B_{\ell,n}u_1,\ldots , p^\ell B_{\ell,n}u_j)+
\wb{B}_{\wb{\ve}}(0)\,,
\end{equation}
where $B_{\ell,n}:=pg'(g^{n-\ell-1}(p^nx))\circ\cdots
\circ pg'(g^{n-N}(p^nx))\in \cL(E)$.
By the same argument,
\[
g^x_s(p^\ell A_{\ell,n}u_1,\ldots,p^\ell A_{\ell,n}u_j)
\in g^x_s(p^\ell A_{\ell,N}u_1,\ldots,A_{\ell,N}u_j)+ \wb{B}_{\wb{\ve}}(0)
\]
for all $j,\ell,s$, and $n$ as before.
Since $g^{n-\nu}(p^nx)\to f(p^\nu x)$
as $n\to\infty$
for all $\nu\in \N$,
we see that $B_{\ell,n}\to A_{\ell,N}$
as $n\to\infty$ and $g^x_{n,s}\to g^x_s$,
for all $\ell$ and $s$ as before.
We therefore
find $N_5\geq N$ such that,
for all $n\geq N_5$, we have
\begin{equation}\label{lle}
\|g^x_s-\wt{g}^x_{n,s}\|\leq \wb{\ve}p^{-j\ell}
\end{equation}
and $\|B_{\ell,n}-A_{\ell,N}\|\leq \wb{\ve}$
for all $j\in \{2,\ldots, k\}$,
$\ell\in \{1,\ldots, N\}$,
and $s\in S_{\ell,j}$.
Then
\[
g^x_s(p^\ell A_{\ell,N}u_1,\ldots,p^\ell A_{\ell,N}u_j)-
\wt{g}^x_{n,s}(p^\ell B_{\ell,n}u_1,\ldots, p^\ell B_{\ell,n}u_j)
\in \wb{B}_{\wb{\ve}}(0)
\]
for all $n\geq N_5$, by (\ref{lle}) and repeated application 
of Lemma~\ref{multilin}.
Combining this with (\ref{end-1}),
we see that (\ref{end-2}) holds
and hence also (\ref{creach}), (\ref{finish}), and (\ref{vani}).
This completes the proof
of Lemma~\ref{mischief} and thus also the proof
of Proposition~\ref{difficult}.
\end{proof}
\section{$\!\!\!$Existence and compatibility
of an analytic~structure}
In this section,
we complete the programme
sketched in the Introduction.
We prove that every finite-dimensional $p$-adic
$C^k$-Lie group
satisfies the hypotheses
of Lazard's Theorem, and show that
Lazard's analytic structure
is $C^k$-compatible with the given $C^k$-manifold structure.
In a final step, we then pass from
$p$-adic Lie groups to the case
of Lie groups over finite extension fields
of $\Q_p$.
\begin{prop}\label{secondkind}
If $\K=\Q_p$
in the situation of Proposition~{\rm \ref{powers}},
$G$ is finite-dimensional,
and the given norm on $L(G)$ is
a maximum norm with respect
to some basis $e_1,\ldots, e_d$ of $L(G)$,
then we can achieve that,
in addition to {\rm (a)--(r)},
also the following holds:
\begin{itemize}
\item[\rm (s)]
For any $s=p^{-j}\leq r$,
the map
\[
\psi\!: (\Z_p)^d\to W_{p^{-1}s}=V_s,\quad
\psi(z_1,\ldots, z_d):=\exp(z_1p^{j+1}e_1)*\cdots*\exp(z_dp^{j+1}e_d)
\]
is a $C^k_{\Q_p}$-diffeomorphism.
\item[\rm (t)]
For any $s\in \;]0,r]$, the open subgroup $\phi^{-1}(V_s)\isom V_s$
of $G$ satisfies
Lazard's conditions {\bf L1--L3},
whence $G$ can be given a finite-dimensional
$p$-adic analytic manifold structure
compatible with the given topology and making it
a $p$-adic analytic Lie group $\wt{G}$.
\item[\rm (u)]
The $p$-adic analytic
manifold structure on $\wt{G}$ is
$C^k_{\Q_p}$-compatible with the $C^k_{\Q_p}$-manifold
structure on~$G$, i.e., $\id\!: \wt{G}\to G$
is a $C^k_{\Q_p}$-diffeomorphism.
\end{itemize}
\end{prop}
\begin{proof}
(s) By Proposition~\ref{filtration}\,(h),
we have $r=p^{-j_0}$ for some $j_0\in \N$.
The map
$f\!: (\Z_p)^d \to W_{p^{-1}r}=V_r$,
$f(z_1,\ldots, z_d):=\exp(z_1p^{j_0+1}e_1)\cdots \exp(z_dp^{j_0+1}e_d)$
is $C^k_{\Q_p}$
by Corollary~\ref{crexp}.
Since $f'(0)(z_1,\ldots,z_d)=p^{j_0+1}(z_1e_1+\cdots+z_de_d)$,
the linear map $p^{-j_0-1}f'(0)$ is a bijective isometry.
Using the Inverse Function Theorem
(\cite{IMP} Prop.\,7.14 and Thm.\,7.3),
we find $i\in \N_0$ such that $f$ induces a $C^k_{\Q_p}$-diffeomorphism
from $p^i\Z_p^d=B_{p^{-i+1}}^{\Z_p^d}(0)$
onto an open subset of $L(G)$,
and such that
$p^{-j_0-1}f$ is an isometry from
the latter set onto $B_{p^{-i+1}}^{L(G)}(0)$.
After replacing $r$ with $p^{-i}r$,
the assertion then holds.

(t) Let $s\in \;]0,r]$; since $V_s=V_{p^{-j}}$
with $j$ chosen such that $p^{-j-1}<s\leq p^{-j}$,
we may assume without loss of generality
that $s=p^{-j}$ for some $j\in \N$.
By Proposition~\ref{filtration} (h)--(j),
$V_s$ is a pro-$p$-group such that $V^{\{p^2\}}=V_{p^{-2}s}\supseteq
[V_s,V_s]$.
With notation as in~(s), we have
$V_s=\psi(\Z_p^d)=\wb{\langle x_1\rangle}\cdots
\wb{\langle x_d\rangle}$, where $x_\nu:= \exp(p^{j+1}e_\nu)$
for $\nu=1,\ldots,d$.
Hence $V_s=\wb{\langle x_1,\ldots,x_d\rangle}$
is finitely generated topologically. Thus conditions
{\bf L1--L3\/} (formulated in the Introduction)
are satisfied by~$V_s$. Then also $U=\psi^{-1}(V_r)$
satisfies these conditions, and thus $G$
can be made a $p$-adic analytic Lie group $\wt{G}$.
The proof of (u)
requires further preparation;
we give it later in this section.
\end{proof}
\begin{defn}\label{defnexp}
Let $G$ be a $C^1$-Lie group
modelled on a topological $\Q_p$-vector space.
A map $\psi\!: \Omega\to G$
on a balanced $0$-neighbourhood
$\Omega\sub L(G)$
(viz.\ $\Z_p\Omega=\Omega$)
is an {\em exponential map\/}
if $\psi(\Omega)$ is an identity neighbourhood,
$\psi$ is continuous at~$0$,
and
$\zeta_x\!: \Z_p\to G$, $\zeta_x(z):=\psi(zx)$
is a homomorphism of class $C^1_{\Q_p}$
with $\zeta_x'(0)=x$,
for each $x\in \Omega$.
\end{defn}
\begin{rem}\label{remsexp}
(a)
Note that, if $\psi\!: \Omega\to G$ is an exponential
map and $x\in \Omega$ with $\psi(x)=1$, then also
$\psi(nx)=1$ for each $n\in \Z$
and hence $\psi(zx)=1$ for each $z\in \Z_p$, by continuity,
entailing that $x=\frac{d}{dt}\big|_{t=0}\psi(tx)=0$.
Hence $\psi$ is injective
and can therefore be considered
as a bijection from $\Omega$ onto the identity neighbourhood
$\psi(\Omega)\sub G$.\vspace{1mm}

(b) If $\psi_i\!: \Omega_i\to G$
are exponential maps for $i\in \{1,2\}$,
then $W:=\psi_1(\Omega_1)\cap\psi_2(\Omega_2)$
is an identity neighbourhood in~$G$.
Given $w\in W$, there exist elements $x_i\in \Omega_i$ such that
$\psi_i(x_i)=w$. Then $\psi_1(nx_1)=\psi_1(x_1)^n=
\psi_2(x_2)^n=\psi_2(nx_2)$
for each $n\in \Z$ and thus
$\psi_1(tx_1)=\psi_2(tx_2)$ for all $t\in \Z_p$,
by continuity, entailing that
$x_1=\frac{d}{dt}\big|_{t=0}\psi_1(tx_1)=
\frac{d}{dt}\big|_{t=0}\psi_2(tx_2)
=x_2$.
Thus $\Omega:=\psi_1^{-1}(W)=\psi_2^{-1}(W)$,
and $\psi_1|_\Omega=\psi_2|_\Omega$. Since $\psi_1$ is continuous,
$\Omega$ is a $0$-neighbourhood. Given $x\in \Omega$,
we have $\psi(x_1)=\psi(x_2)$ and thus
$\psi_1(zx)=\psi_2(zx)$ for each $z\in \Z_p$,
entailing that $\Z_p x\sub \Omega$. Thus $\Omega$ is balanced.

(c) If $G$ is a $p$-adic $C^k$-Lie group
admitting an exponential map $\exp_G\!: \Omega\to G$,
we let
$\Gamma(G)$ be the set of germs $[\xi]$
at~$0$
of continuous homomorphisms $\xi\!:W\to G$
defined on some open subgroup $W\sub \Q_p$.
Then
\begin{equation}\label{thet}
\theta_G\!: L(G)\to \Gamma(G),\quad
\theta_G(x)=[t\mto \exp_G(tx)]
\end{equation}
is a bijection
$[$Since $x=\frac{d}{dt}\big|_{t=0}\exp_G(tx)$
for $x\in L(G)$,
the map $\theta_G$ is injective.
To see that $\theta_G$ is surjective,
let $[\xi]\in \Gamma(G)$,
where $\xi\!: W\to G$.
Since $\xi$ is continuous,
we find $n\in\N$ such that
$p^n\in W$ and $\xi(p^n)\in \exp_G(\Omega)$.
Hence $\xi(p^n)=\exp_G(x)$ for some $x\in \Omega$,
whence $\xi(mp^n)=\exp_G(mx)$ for all
$m\in \Z$ and actually for all $m\in \Z_p$,
by continuity.
Thus $\xi(z)=\exp_G(zp^{-n}x)$ for all $z$
in the $0$-neighbourhood $p^n\Z_p$ in $\Q_p$,
whence $[\xi]=\theta_G(p^{-n}x)$.\,]
We give $\Gamma(G)$ the $p$-adic topological
vector space structure making
$\theta_G$ an isomorphism of topological vector spaces.
Apparently, the map $\theta_G$ only depends on the germ of $\exp_G$ at~$0$
and therefore is independent of the choice
of exponential map, by (b).
\end{rem}
If $G$ is an ultrametric $p$-adic Banach-Lie group
of class $C^{k_+}_{\Q_p}$ and
$\phi\!: U\to V_r\sub L(G)$
and $\exp\!: V_r\to V_r$ are as in Proposition~\ref{powers},
then clearly $\exp_G\!: V_r\to G$, $x\mto \phi^{-1}(\exp(x))$
is an exponential map for~$G$.
We define $\log_G:=(\exp_G|^U)^{-1}=
\log\circ \phi^{-1}\!: U\to V_r\sub L(G)$.
Given $x\in U^{\{p^n\}}=\phi^{-1}(V_{p^{-n}r})$
(see Proposition~\ref{filtration}\,(i)),
we let $x^{p^{-n}}:=\phi^{-1}(\tau_p^{-n}(\phi(x)))$;
this is the unique element in $U$ with $p^n$-th power~$x$.
Using a version of
Trotter's Product Formula,
we now show that the topological vector space $\Gamma(G)$
is determined by the topological group underlying~$G$.
\begin{la}\label{trotter}
In the preceding situation, we have:
\begin{itemize}
\item[\rm (a)]
$x+y=\lim_{n\to \infty}
\log_G\big( (\exp_G(p^nx)\exp_G(p^ny))^{p^{-n}}\big)$
for all $x,y\in V_r$.
\item[\rm (b)]
If $\,[\xi_1]$, $[\xi_2]\in \Gamma(G)$
and $[\xi_1]+[\xi_2]=[\xi_3]$,
there is $\delta>0$
such that $\xi_1(t)$, $\xi_2(t)$ and $\xi_3(t)$
are defined
for all $\,t\in \Z_p$ with $|t|_p\leq \delta$, each of them
is an element of $U$, and
$\xi_3(t)=\lim_{n\to\infty}(\xi_1(p^nt)\xi_2(p^nt))^{p^{-n}}$.
\item[\rm (c)]
If $\, \wt{G}$ is a $p$-adic ultrametric
Banach-Lie group of class $C^{k_+}_{\Q_p}$
such that $\wt{G}$ and $G$ have
the same underlying topological group,
then $\Gamma(G)=\Gamma(\wt{G})$
as a set and as a topological vector
space over~$\Q_p$.
Furthermore, the map
$G\to\wt{G}$, $x\mto x$
is an isomorphism of $SC^1_{\Q_p}$-Lie groups.
\end{itemize}
\end{la}
\begin{proof}
(a) Since $G\isom V_r$, we may assume
that $G=V_r$.
Note that, for $v\in V_{p^{-n}r}$,
we have $\exp(p^n(\log(v^{p^{-n}})))
=(\exp(\log(v^{p^{-n}})))^{p^n}=(v^{p^{-n}})^{p^n}=v=\exp(\log(v))$,
whence $p^n\log(v^{p^{-n}})=\log(v)$ and
thus $\,\log(v^{p^{-n}})=p^{-n}\log(v)$.
Therefore
\begin{equation}\label{entay}
\log((\exp(p^nx)\exp(p^ny))^{p^{-n}})=p^{-n}\log(\exp(p^nx)\exp(p^ny))\quad
\mbox{for all $x,y\in V_r$, $n\in \N$.}
\end{equation}
Now $f\!: \Z_p\to L(G)$, $f(z):=\log(\exp(zx)\exp(zy))$
being $C^1$ and thus strictly differentiable
at~$0$,
with $f'(0)=x+y$, given $\ve>0$
we find $\delta>0$ such that
$\|f(z)-z(x+y)\|\leq \ve |z|_p$
for all $z\in \Z_p$ such that $|z|_p\leq\delta$.
Using (\ref{entay}), we obtain
$\|\log((\exp(p^nx)\exp(p^ny))^{p^{-n}})-(x-y)\|=
p^n \|f(p^n)-p^n(x+y)\|\leq \ve$ for all $n\in \N$ such that
$p^{-n}\leq\delta$. Hence (a) holds.

(b) Let $v_i:=\theta_G^{-1}([\xi_i])$
for $i\in \{1,2,3\}$.
Then $v_3=v_1+v_2$, the map $\theta_G$
being linear.
There is $\delta>0$
such that
$\delta\|v_i\|<r$
for each $i\in \{1,2,3\}$,
$\xi_i(t)$ is defined for
all $t\in \Z_p$ with $|t|_p\leq \delta$,
and $\xi_i(t)=\exp_G(tv_i)$.
Using Part\,(a) and
the continuity of $\exp_G$,
for any $t$ as before we obtain
$\xi_3(t)=\exp_G(tv_3)=\exp_G(tv_1+tv_2)
=\lim_{n\to\infty}(\exp(p^ntv_1)\exp(p^ntv_2))^{p^{-n}}
=\lim_{n\to\infty}(\xi_1(p^nt)\xi_2(p^nt))^{p^{-n}}$,
as asserted.

(c) Because the definition of the set
$\Gamma(G)$ only involves the topological group structure
of~$G$, we have $\Gamma(G)=\Gamma(\wt{G})$ as a set.
Given $z\in \Q_p$ and $[\xi]\in \Gamma(G)$,
we have $z[\xi]=[t\mto \xi(zt)]$
both in $\Gamma(G)$ and $\Gamma(\wt{G})$.
Thus the scalar multiplication maps of the two
$\Q_p$-vector
spaces coincide.
We let $\wt{U}=\wt{\phi}^{-1}(\wt{V}_{\wt{r}})\sub \wt{G}$
be an open subgroup
of~$\wt{G}$ playing a role analogous to that
of $U\sub G$; after shrinking $\wt{r}$,
we may assume without loss of generality that
$\wt{U}\sub U$.
Then $\wt{U}^{\{p^n\}}$
is an open subgroup of $G$ for each
$n\in \N$, and each $x\in \wt{U}^{\{p^n\}}$
has a unique $p^n$-th root $y$
in $\wt{U}$,
and also a unique $p^n$-th root $x^{p^{-n}}$
in $U$; by uniqueness of $x^{p^{-n}}$ in $U$,
we have $x^{p^{-n}}=y\in \wt{U}$.
Now assume that $[\xi_1],[\xi_2]\in \Gamma(G)$
are given, with
$[\xi_1]+[\xi_2]=[\xi_3]$ in $\Gamma(G)$,
$[\xi_1]+[\xi_2]=[\xi_4]$ in $\Gamma(\wt{G})$.
Part\,(b) allows us to calculate the respective sum
both in
$\Gamma(G)$
and $\Gamma(\wt{G})$:
there is $\delta>0$ such that
$\xi_1(t)$, $\xi_2(t)$, $\xi_3(t)$ and $\xi_4(t)$
are defined for all $t\in \Q_p$ with $|t|_p\leq \delta$,
all of them are elements of $\wt{U}$
(hence of $U$),
and
\[
\xi_3(t)=\lim_{n\to\infty}(\xi_1(p^nt)\xi_2(p^nt))^{p^{-n}}
=\xi_4(t)\,,
\]
using that, as just explained,
the $p^n$-th roots in $U$ and $\wt{U}$
occurring here coincide.
By the preceding, the germs at $0$ of $\xi_3$ and $\xi_4$
coincide. Thus $[\xi_3]=[\xi_4]$,
whence the sum $[\xi_1]+[\xi_2]$ is the same
in $\Gamma(G)$ and $\Gamma(\wt{G})$.
Hence $\Gamma(G)$ and $\Gamma(\wt{G})$
coincide as $\Q_p$-vector spaces.

To see that $\Gamma(G)$ and $\Gamma(\wt{G})$
coincide as topological vector spaces,
consider
the open $0$-neighbourhoods
$\Omega:=\theta_G(V_r)\sub \Gamma(G)$
and $\wt{\Omega}:=\theta_{\wt{G}}(\wt{V}_{\wt{r}})\sub \Gamma(\wt{G})$.
The maps
$\Log_G:=\theta_G|_{V_r}^\Omega \circ \log_G\!:U\to \Omega$
and
$\Log_{\wt{G}}:=\theta_{\wt{G}}|_{\wt{V}_{\wt{r}}}^{\wt{\Omega}}\circ
\log_{\wt{G}}\!:\wt{U}\to \wt{\Omega}$
are $SC^1$-diffeomorphisms.
For each $x\in U\cap \wt{U}$,
there is a unique continuous homomorphism $\gamma_x\!:
\Z_p\to G$ such that $\gamma_x(1)=x$
(namely, $\gamma_x=\phi^{-1}\circ \eta_{\phi(x)}$).
Since $\exp_G(n\log_G(x))=\exp_G(\log_G(x))^n=x^n$
for each $n\in \Z$, we deduce that
$\gamma_x(z)=\exp_G(z\log_G(x))$
for all $z\in \Z_p$, and likewise
$\gamma_x(z)=\exp_{\wt{G}}(z\log_{\wt{G}}(x))$.
As a consequence,
\[
\Log_G(x)=\theta_G(\log_G(x))=[t\mto \exp_G(t\log_G(x))]=[\gamma_x]
=\theta_{\wt{G}}(\log_{\wt{G}}(x))=\Log_{\wt{G}}(x)\,,
\]
entailing that
$x\mto \Log_{\wt{G}}(\Log_G^{-1}(x))=x$
is an $SC^1$-diffeomorphism (and hence a homeomorphism)
from the open $0$-neighbourhood
$Q:=\Log_G(U\cap \wt{U})$ in $\Gamma(G)$
onto the open identity neighbourhood $\Log_{\wt{G}}(U\cap\wt{U})$
in $\Gamma(\wt{G})$. As a consequence,
$\Gamma(G)=\Gamma(\wt{G})$ as a topological $\Q_p$-vector space.
Since $\Log_G|_{U\cap \wt{U}}^Q=\Log_{\wt{G}}|_{U\cap \wt{U}}^Q$
is an $SC^1$-diffeomorphism
both on $U\cap \wt{U}$, considered
as an open subset of $G$, and as
an open subset of~$\wt{G}$,
we readily deduce that both the homomorphism
$\id\!: G\to \wt{G}$
and its inverse $\id\!: \wt{G}\to G$ are of class
$SC^1$.
\end{proof}
{\bf Proof of Proposition~\ref{secondkind}, completed.}
(u) By Proposition~\ref{secondkind}\,(t),
there is a
$p$-adic analytic manifold structure on $G$,
compatible with the given topology, which makes
$G$ a finite-dimensional, $p$-adic analytic Lie group
$\wt{G}$. By Lemma~\ref{trotter}\,(c),
the analytic Lie group structure on $\wt{G}$ is $C^1$-compatible
with the given $C^k$-manifold structure
on~$G$.
By Part\,(s) of Proposition~\ref{secondkind},
the map
\[
\psi\!: (\Z_p)^d\to V_r,\quad
\psi(z_1,\ldots,z_d):=\zeta_1(z_1)\zeta_2(z_2)\cdots \zeta_d(z_d)
\]
is a $C^k$-diffeomorphism, where
$\zeta_\nu\!: \Z_p\to G$,
$\zeta_\nu:=\exp_G(zp^{j_0+1}e_\nu)$ for $\nu\in \{1,\ldots, d\}$.
Each $\zeta_\nu$ is, in particular, a continuous
homomorphism and hence analytic
as a map into $\wt{G}$
by Cartan's Theorem
(\cite{Ser}, Part~II, Chapter~V,
\S9, Thm.\,2). Hence $\psi$ is analytic as
a map into $\wt{G}$.
Now, $\psi$ being a $C^1$-diffeomorphism onto
$U$ considered as
an open subset of~$G$,
the map $\psi$ also is a $C^1$-diffeomorphism onto
$U$ considered as an open subset of $\wt{G}$,
the two manifold structures
being $C^1$-compatible.
Being a $C^1$-diffeomorphism and analytic,
$\psi\!: (\Z_p)^d\to U\sub \wt{G}$
is an analytic diffeomorphism
(and hence a $C^k$-diffeomorphism),
as a consequence of the
Inverse Function Theorem for
analytic maps \cite[p.\,73]{Ser}.
Thus both $G$ and $\wt{G}$ induce the same $C^k$-manifold
structure on the open identity neighbourhood~$U$,
whence the homomorphisms $\id\!: G\to \wt{G}$
and $\id\!: \wt{G}\to G$ are $C^k$,
being $C^k$ on $U$
(Lemma~\ref{higherdhom}).\vspace{2.5mm}\Punkt

\noindent
We now prove
our main result, Theorem~A (from the Introduction):\\[2.5mm]
{\bf Proof of Theorem~A.}
Being a $C^k_\K$-Lie group, $G$ can also be considered
as a $C^k_{\Q_p}$-Lie group.
Thus Proposition~\ref{secondkind}
provides a finite-dimensional
$p$-adic analytic manifold structure on $G$
making it a $p$-adic analytic Lie group $\wt{G}$,
which is $C^k_{\Q_p}$-compatible
with the given $C^k_{\Q_p}$-manifold structure
on~$G$. Then $L(\wt{G})=T_1\wt{G}$
can be identified with $L(G)$,
considered as $\Q_p$-vector space,
in a natural way.
Given $x\in G$, consider the inner automorphism
$I_x\!: G\to G$, $I_x(y):=xyx^{-1}$
of the $C^k_\K$-Lie group~$G$
and the corresponding
$\K$-linear tangent map
$\Ad_x:=L(I_x):=T_1(I_x)\!: L(G)\to L(G)$.
Obviously the same mapping $\Ad_x$
is obtained when considering $I_x$
as an automorphism of $\wt{G}$,
and so the given $\K$-vector space structure
on $L(G)$ is compatible with the adjoint action of
$\wt{G}$. Furthermore, the Lie bracket on $L(G)$
as the Lie algebra of $\wt{G}$
is $\K$-bilinear with respect to the
given $\K$-vector space structure
on $L(G)$. To see this, note
that the image of
the $\Q_p$-analytic homomorphism
$h\!: \wt{G}\to \GL_{\Q_p}(L(G))$,
$h(x):=\Ad_x$
is contained in the closed subgroup
$\GL_\K(L(G))$,
whence the image of
$\ad:=L(h)\!:
L(\wt{G})\to \gl_{\Q_p}(L(G))$
is contained in the corresponding Lie subalgebra
$\gl_\K(L(G))$ of $\K$-linear endomorphisms.
Thus $\ad(x).y=[x,y]$ is $\K$-linear in $y$ for each
$x\in L(G)$, and hence so it is in $x$,
by antisymmetry of the Lie bracket.
Applying \cite{BLi},
Chapter~III, \S4.2, Cor.\,2 to Thm.\,2,
we now obtain a unique $\K$-analytic manifold structure
on $G$ making it a $\K$-analytic Lie group
$\wh{G}$ with $\K$-Lie algebra $L(G)$,
and $\Q_p$-analytically compatible
with the $p$-adic analytic structure
on $\wt{G}$.
We let $\exp_{\wh{G}}\!: \Omega\to \wh{G}$
be a $\K$-analytic exponential map for
the $\K$-analytic Lie group $\wh{G}$,
in the sense of \cite{BLi}, Chapter~III, \S4.3,
Definition~1, defined on an open $\A$-submodule $\Omega$
of $L(G)$ (where, as before, $\A=\{z\in \K\!: |z|\leq 1\}$).
By {\em loc.cit.},
Thm.\,4\,(i),
there exists $\ve>0$ such that
$\exp_{\wh{G}}((t+s)x)=\exp_{\wh{G}}(tx)\exp_{\wh{G}}(sx)$
for all $x\in \Omega$ and all $s,t\in \K$ such that
$|s|,|t|\leq \ve$.
After replacing $\Omega$ by $p^n\Omega$
with $n$ sufficiently large, we may assume that $\ve=1$ here,
entailing that $\exp_{\wh{G}}(zx)=(\exp_{\wh{G}}(x))^z$
for all $z\in \Z$ and hence also for all
$z\in \Z_p$. Since $\exp_{\wh{G}}'(0)=\id_{L(G)}$,
after shrinking $\Omega$ we may assume that $\exp_{\wh{G}}$
is a $\K$-analytic diffeomorphism onto an open subset of~$\wh{G}$.
Note that, by the last and penultimate property,
$\exp_{\wh{G}}$ also is an exponential
map in the sense of Definition~\ref{defnexp}
for the $C^1_{\Q_p}$-Lie group underlying $\wh{G}$,
and hence for $G$.
On the other hand, Proposition~\ref{powers}\,(r)
provides an exponential map $\exp_G:=\phi^{-1}\circ \exp\!:
V_r\to U$ for $G$, such that\vspace{-4mm}
\[
\zeta_x\!: \A\to G,\quad z\mto \exp_G(zx)
\]
is of class $C^k_\K$, for each $x\in V_r$
(Corollary~\ref{crexp}).
Since both $\exp_{\wh{G}}$ and $\exp_G$
are exponential maps for $G$, considered as an $C^1_{\Q_p}$-Lie group,
we deduce from Remark~\ref{remsexp}\,(b)
that $\exp_G$ and $\exp_{\wh{G}}$
coincide on $Q:=B_s^{L(G)}(0)\sub \Omega$
for some $s\in \;]0,r]$.
As a consequence, for every $x\in Q$
the map $\zeta_x=\exp_G(\sbull x)=\exp_{\wh{G}}(\sbull x)$
is $C^k_\K$ both as a map into $G$, and as a map
into $\wh{G}$.
We pick a basis $e_1,\ldots, e_d\in Q$
of the $\K$-vector space $L(G)$.
Using the Inverse Function Theorem for
$C^k_\K$-maps and the Inverse Function Theorem
for $\K$-analytic maps, we find $n\in \N_0$
such that
\[
(p^{-n}\A)^d\to G,\quad (z_1,\ldots,z_d)\mto
\zeta_{e_1}(z_1)\cdots \zeta_{e_d}(z_d)
\]
is both a $C^k_\K$-diffeomorphism onto an open subset
of $G$, and a $\K$-analytic diffeomorphism
onto the corresponding subset of $\wh{G}$.
Hence $\id\!: G\to \wh{G}$ is $C^k_\K$
on some open identity neighbourhood
and thus $C^k_\K$, and likewise for $\id\!:
\wh{G}\to G$.\Punkt
\appendix
\section{Proofs of the lemmas from Section~1}
{\bf Proof of Lemma~\ref{taylor}.}
The map $f$ being $C^2$, we have a second order Taylor expansion
\[
f(x+ty)-f(x)-tdf(x,y)=t^2 a_2(x,y)+t^2\,R_2(x,y,t)
\]
for $(x,y,t)\in U^{[1]}$,
with remainder $R_2\!: U^{[1]}\to F$
(see \cite{BGN}, Thm.\,5.1 and Prop.\,5.3).
Let $\|.\|_\gamma$ be a continuous seminorm on~$F$.
Since $R_2(x_0,0,0)=0$ and $a_2(x_0,0)=0$,
there exists $\rho\in \;]0,1]$
such that $B_{2\rho}(x_0)\sub U$,
\[
\|R_2(x,y,t)\|_\gamma \leq 1\qquad
\mbox{for all $x\in B_\rho(x_0)$, $y\in B_\rho(0)$,
and $|t|<\rho$,}
\]
and
$\|a_2(x,y)\|_\gamma \leq 1$ for all $x\in B_\rho(x_0)$ and
$y\in B_\rho(0)$.
Pick $a\in \K^\times$ such that $|a|<1$;
define $\delta:=\rho^2|a|<\rho$ and $C:=2/(\rho|a|)^2$.
Let $x\in B_\delta(x_0)$ and $y\in B_\delta(0)$.
If $y=0$, then $\|f(x+y)-f(x)-df(x,y)\|_\gamma \leq C\|y\|^2$
trivially. If $y\not=0$, there exists
$k\in \Z$ such that
$|a|^{k+1}\leq \rho^{-1}\|y\| <|a|^k$.
Then $\|a^{-k}y\|<\rho$ and $|a^k|\leq |a|^{-1}\rho^{-1}\|y\|<\rho$, and thus
\[
f(x+y)-f(x)=f(x+a^ka^{-k}y)-f(x)=df(x,y)
+a^{2k}a_2(x,a^{-k}y)+a^{2k}R_2(x,a^{-k}y,a^k)
\]
where
$\|a^{2k}a_2(x,a^{-k}y)+a^{2k}R_2(x,a^{-k}y,a^k)\|_\gamma \leq |a|^{2k}
(\|a_2(x,a^{-k}y)\|_\gamma
+\|R_2(x,a^{-k}y,a^k)\|_\gamma )\leq 2|a|^{2k}\leq 2|a|^{-2}\rho^{-2}\|y\|^2
=C\|y\|^2$.\vspace{2mm}\Punkt

\noindent
{\bf Proof of Lemma~\ref{taylor2}.}
We may assume that $U$ and $V$ are balanced.
Then, for every $x\in U$, $y\in V$ and $s,t\in \K$
such that $|s|,|t|\leq 1$, the first order Taylor
expansion of $f$ about $(tx,0)$ shows that
\begin{equation}\label{expa1}
f(tx,sy)=\underbrace{f(tx,0)}_{=f(0,0)+\lambda(tx)}+sdf((tx,0),(0,y))+s
R_1((tx,0),(0,y),s)\,,
\end{equation}
where $df\!: U\times V\times E\times F\to H$
and the remainder
$R_1\!: (U\times V)^{[1]}\to H$
are of class $C^1_\K$ (see \cite{BGN}).
Using the first order Taylor expansions of $df$
and $R_1$ about $((0,0),(0,y))$ and $((0,0),(0,y),s)$,
with remainders $P_1$ and $Q_1$,
respectively, we deduce from (\ref{expa1}) that
\begin{eqnarray}
f(tx,sy) &=& f(0,0)+\lambda(tx)+
s\,\underbrace{df((0,0),(0,y))}_{=\mu(y)}
+st\, \underbrace{d(df)((0,0,0,y),(x,0,0,0))}_{=d^2f((0,0),(0,y),(x,0))}\nonumber \\
& &\;\; +\; st\, P_1((0,0,0,y),(x,0,0,0),t) +
s\,\underbrace{R_1((0,0),(0,y),s)}_{=0}\nonumber \\
& &\;\; +\; st\, dR_1((0,0,0,y,s),(x,0,0,0,0))
+st\, Q_1((0,0,0,y,s),(x,0,0,0,0),t)\nonumber \\
&=& f(0,0)+\lambda(tx)+\mu(sy)+\beta(tx,sy)+ st\, g(x,y,s,t)\,,\label{betg}
\end{eqnarray}
where inessential brackets were
suppressed in the notation,
the 6th term vanishes due to (\ref{cc2}),
the map $\beta\!: E\times F\to H$, $\beta(u,v):=d^2f((0,0),(0,v),(u,0))$
is continuous bilinear,
and where the 5th, 7th and 8th terms are combined
in an apparent way in the form $st\, g(x,y,s,t)$,
where $g\!: U\times V\times \wb{B}_1^\K(0)\times \wb{B}_1^\K(0)\to H$
is continuous, and $g(0,0,0,0)=0$.
We find $\sigma\in \,]0,1]$ such that $B_\sigma^E(0)\sub U$,
$B_\sigma^F(0)\sub V$, and $\|g(x,y,s,t)\|\leq 1$
for all $x\in B_\sigma^E(0)$, $y\in B_\sigma^F(0)$,
and $s,t\in \K$ such that $|s|,|t|\leq \sigma$.
Let $a\in \K^\times$ be an element such that
$|a|<1$. Assume that
$x\in E$ and $y\in F$ such that
$\|x\|, \|y\|< \sigma^2|a|=:\delta$.
If $x=0$ or $y=0$, then $\|f(x,y)-f(0,0)-\lambda(x)-\mu(y)\|_\gamma=0$
by (\ref{cc1}) and (\ref{cc2}), respectively.
Otherwise, we find uniquely determined
numbers $k,\ell\in \Z$
such that $|a|^{k+1}\leq \frac{\|x\|}{\sigma} <|a|^k$
and $|a|^{\ell+1}\leq \frac{\|y\|}{\sigma} < |a|^\ell$.
Then $|a|^k \leq \frac{\|x\|}{|a|\sigma}<\sigma$
and $\|a^{-k}x\|<\sigma$;
similarly, $|a|^\ell<\sigma$ and $\|a^{-k}y\|<\sigma$.
Choosing $t:=a^k$ and $s:=a^\ell$ in (\ref{betg}),
we obtain
\begin{eqnarray*}
\|f(x,y)-f(0,0)-\lambda(x)-\mu(y)\|_\gamma & =& \|f(a^k(a^{-k}x),a^\ell(a^{-\ell}y))-f(0,0)-\lambda(x)-\mu(y)\|_\gamma\\
&\leq &\|\beta(x,y)\|_\gamma +
|a|^k|a|^\ell\|g(a^{-k}x, a^{-\ell}y,a^\ell,a^k)\|_\gamma\\
& \leq &
\|\beta\|_\gamma
\|x\|\,\|y\| + \frac{\|x\|}{|a|\sigma}\,\frac{\|y\|}{|a|\sigma}
\leq C\,\|x\|\,\|y\|,
\end{eqnarray*}
with $C:=\|\beta\|_\gamma + (|a|\sigma)^{-2}$.\vspace{3mm}\Punkt

\noindent
The proof of Lemma~\ref{smallo} will be based on
the following observation:
\begin{la}\label{subla}
Let $X$ and $E$ be normed spaces
over a valued field~$\K$, $\,F$ be a polynormed
$\K$-vector space,
$U\sub X$ an open
subset, $n\in \N$, and $f\!: U\times E^n\to F$
be a $C^1_\K$-map
such that $f(x,\sbull)\!: E^n\to F$ is $n$-linear,
for each $x\in U$. Then the map
$\phi\!: U\to \cL^n(E,F)$, $\phi(x):=f(x,\sbull)$
is continuous.
\end{la}
\begin{proof}
Using the first order Taylor expansion
of $f$, we can write
\[
f(x+tz,y)-f(x,y)=td_1f(x,y,z)
+tR(x,z,t,y)\quad\mbox{for all $(x,z,t)\in U^{[1]}$
and $y\in E^n$,}
\]
where $d_1f(x,y,z):=df((x,y),\, (z,0))$
is linear in $z$
and $R\!: U^{[1]}\times E^n\to F$
is a continuous map such that $R(x,z,t,y)=0$
whenever $t=0$. Since $f(x,\sbull)$ is $n$-linear,
apparently so is $d_1f(x,\sbull,z)$,
and then also $R(x,z,t,\sbull)$
(this is clear for $t\not=0$,
and follows for $t=0$ by continuity).
Thus $d_1f(x,\sbull)$ is a continuous
$(n+1)$-linear map.
Now assume that $x\in U$, $\ve>0$,
and that $\|.\|_\gamma$ is a continuous seminorm on~$F$.
Then there is $\delta\in \;]0,1]$ such that $B_\delta^X(x)\sub U$
and
$\|R(x,z,t,y_1,\ldots,y_n)\|_\gamma\leq 1$
for all $z\in B_\delta^X(0)$, $|t|\leq \delta$, and
$y_1,\ldots, y_n\in B_\delta^E(0)$.
Pick $a\in \K^\times$ such that $|a|<1$.
Let $\rho:=\min\{\ve (1+2\|d_1f(x,\sbull)\|_\gamma)^{-1},\frac{\ve}{2}
(|a|\delta)^{n+1}, |a|\delta^2\}$.
Let $0\not=z\in B_\rho^X(0)$
and $y_1,\ldots,y_n\in E$.
If some $y_j=0$, then $(\phi(x+z)-\phi(x))(y_1,\ldots,y_n)=0$.
Otherwise, we find $k_0,k_1,\ldots, k_n\in \Z$
such that $|a|^{k_0+1}\leq \|z\|\delta^{-1}<|a|^{k_0}$
and $|a|^{k_j+1}\leq \|y_j\|\delta^{-1}<|a|^{k_j}$
for $j=1,\ldots, n$.
Thus $\|a^{-k_0}z\|<\delta$ and $\|a^{-k_j}y_j\|<\delta$
and hence, writing $y:=(y_1,\ldots, y_n)$,
\begin{eqnarray}
f(x+z,y)-f(x,y)\!\!&\!=\!& \!\!f(x+a^{k_0} a^{-k_0}z,y)-f(x,y)
=d_1f(x,y,z)+a^{k_0}R(x,a^{-k_0}z,a^{k_0},y)\nonumber \\
\!\!&\!=\!&\!\!d_1f(x,y,z)+a^{k_0+k_1+\cdots +k_n}R(x,a^{-k_0}z,a^{k_0},a^{-k_1}y_1,\ldots,
a^{-k_n}y_n),\label{secline}
\end{eqnarray}
where
$\|d_1f(x,y,z)\|_\gamma\leq \|d_1f(x,\sbull)\|_\gamma\cdot \|z\|\cdot \|y_1\|\cdots
\|y_n\|\leq \frac{\ve}{2}\|y_1\|\cdots\|y_n\|$.
Since
\[
|a^{k_0+k_1+\cdots+k_n}|\leq
|a|^{-(n+1)}\delta^{-(n+1)}\|z\|\cdot\|y_1\|\cdots\|y_n\|
\leq \frac{\ve}{2}\|y_1\|\cdots\|y_n\|
\]
and $\|R(x,a^{-k_0}z,a^{k_0},a^{-k_1}y_1,\ldots,
a^{-k_n}y_n)\|_\gamma\leq 1$,
also the final term in (\ref{secline})
has norm $\leq \frac{\ve}{2}\|y_1\|\cdots\|y_n\|$.
Thus $\|(\phi(y+z)-\phi(x))(y_1,\ldots,y_n)\|_\gamma
\leq \ve\|y_1\|\cdots\|y_n\|$
for all $y_1,\ldots,y_n\in E$ and hence
$\|\phi(x+z)-\phi(x)\|_\gamma \leq\ve$, for all $z\in B_\rho^X(0)$.
Thus $\phi$ is continuous.
\end{proof}
\noindent
{\bf Proof of Lemma~\ref{smallo}.}
If $f$ is of class $C^{k+1}$,
then $f^{[k]}$ is of class $C^1$
and hence so is $d^kf\!: U\times E^k\to F$,
being a partial map of $f^{[k]}$.
Since $d^kf(x,\sbull)$ is $k$-linear
for each $x\in U$, Lemma~\ref{subla}
shows that $U\to \cL^k(E,F)$, $x\mto d^kf(x,\sbull)$
is continuous.\\[3mm]
If $f$ is merely $C^k$ but $\K$ a complete
valued field and
$E$ finite-dimensional, then
we pick a basis $e_1,\ldots,e_n$
of $E$ and recall that $\K^n\to E$, $(t_1,\ldots,t_n)\mto \sum_{j=1}^n t_je_j$
is an isomorphism of topological vector spaces
(\cite{BTV}, Chapter~I, \S2, No.\,3, Thm.\,3).
The continuity of $\phi$ now follows from the continuity
of the maps $U\to F$, $x\mto
d^kf(x,e_{i_1},\ldots, e_{i_k})$
and the (easily verified) fact
that $\cL^k(E,F)\to F^{\{1,\ldots,n\}^k}$,
$\beta\mto (\beta(e_{i_1},\ldots, e_{i_k}))_{i_1,\ldots,i_k=1}^n$
is an isomorphism of topological vector spaces.\vspace{2mm}\Punkt

\noindent
{\bf Proof of Lemma~\ref{ntayl}.}
Let $z\in U$,
$\|.\|_\gamma$ be a continuous seminorm on~$F$,
and $\ve\in \;]0,1]$.\\[3mm]
{\em If $f$ is of class $C^{k+1}$\/},
then an apparent adaptation of the proof of Lemma~\ref{taylor}
based on the $(k+1)$-th order Taylor expansion of
$f$ gives $\delta\in \;]0,1]$ and $C>0$ such that
$B_{2\delta}(z)\sub U$ and
\[
\left\|f(x+y)-f(x)-
\!\!\sum_{j=1}^k\frac{1}{j!}d^jf(x,y,\ldots,y)\right\|_\gamma
\leq C\|y\|^{k+1}\quad\mbox{for all $x\in B_\delta(z)$
and $y\in B_\delta(0)$.}
\]
Set $\sigma:=\frac{\delta\ve}{1+C}<\delta$.
Then $\|R(x,y)\|_\gamma\leq C\|y-x\|^{k+1}\leq \ve \|y-x\|^k$
for all $x,y\in B_\sigma(z)$. We deduce that (\ref{remaind})
holds.\\[3mm]
{\em If $f$ is $C^k$, $\K$ is locally compact and $E$
finite-dimensional}, let $R_k\!: U^{[1]}\to F$ be the remainder
of the $k$-th order Taylor expansion.
Choose $\rho>0$ such that $\wb{B}_{2\rho}(z)\sub U$.
Pick $a\in \K^\times$ such that $|a|<1$.
Since $R_k(x,y,0)=0$ for all $x\in U$ and $y\in E$,
using the compactness of $\wb{B}_\rho(z)\times \wb{B}_1(0)$
we find $\delta\in \;]0,\rho]$
such that $\|R_k(x,y,t)\|_\gamma\leq
\ve|a|^k$ for all $(x,y,t)\in
\wb{B}_\rho(z)\times \wb{B}_1(0)\times B_\delta(0)
\sub U^{[1]}$.
Given $x,y \in B_{|a|\delta/2}(z)$,
with $x\not=y$ to avoid trivialities,
there exists $\ell\in \Z$ such that
$|a|^{\ell+1}\leq \|y-x\|<|a|^\ell$.
Since $R(x,y)=R(x,x+(y-x))=R(x,x+a^\ell a^{-\ell}(y-x))=a^{k\ell}
R_k(x,a^{-\ell}(y-x),a^\ell)$,
where $|a^\ell|\leq |a|^{-1}\|y-x\|<\delta$ and
$\|a^{-\ell}(y-x)\|<1$, we deduce that
$\|R(x,y)\|_\gamma=|a|^{k\ell}\|R_k(x,a^{-\ell}(y-x),a^\ell)\|_\gamma \leq
|a|^{-k}\|y-x\|^k\ve |a|^k=\ve \|y-x\|^k$.
Hence (\ref{remaind}) holds.\vspace{2mm}\Punkt

\noindent
{\bf Proof of Lemma~\ref{alongcurv}.}
For each $j=1,\ldots, k$,
and $x\in U$, let $b_j(x,\sbull)\!: E^j\to F$ be the symmetric
$j$-linear map associated with $a_j(x,\sbull)$;
since $b_j(x,\sbull)$ can be obtained from $a_j(x,\sbull)$
by polarization (cf.\ \cite{BaS}),
it is easy to see that $b_j\!: U\times E^j\to F$ is continuous.
If $\eta\!: I \to U$ is a $C^k$-curve, defined on an open
subset $I\sub \K$, we have the Taylor expansion
\begin{equation}\label{expeta}
\eta(s)-\eta(t)=\sum_{i=1}^kc_i(t)(s-t)^j
+r(t,s)(s-t)^k\,,
\end{equation}
where $c_i:=\frac{1}{i!}\eta^{(i)}$ is continuous,
and $r\!: I\times I\to E$ is a continuous map
vanishing on the diagonal.
Substituting (\ref{expeta}) into
$f(y)-f(x)=\sum_{j=1}^kb_j(x,y-x,y-x,\ldots,y-x) +R(x,y)$,
we find that\vspace{-1mm}
\begin{equation}\label{vorlage}
f(\eta(s))-f(\eta(t))
=\sum_{\ell=1}^k g_\ell(t)(s-t)^\ell +(s-t)^k\rho(t,s)\,,
\end{equation}
where\vspace{-11mm}
\[
g_\ell(t)=\sum_{j=1}^k
\sum_{\stackrel{i_1,\ldots,i_j\in \{1,\ldots,k\}}{{\rm s.t.}\;
i_1+\cdots+i_j=\ell}}b_j(\eta(t),c_{i_1}(t),\ldots,
c_{i_j}(t))
\]
and $\rho$ is a sum of terms
of the following form:
Firstly, we have summands of the form
\[
h(t,s)=(s-t)^{i_1+\cdots+i_j-k}b_j(\eta(t),c_{i_1}(t),\ldots,
c_{i_j}(t))\,,\]
where $j\in \{1,\ldots, k\}$
and $i_1,\ldots,i_j\in \{1,\ldots, k\}$
such that $i_1+\cdots+i_j>k$;
any such summand is a continuous map $I\times I\to F$,
and vanishes on the diagonal.
Second, we have summands obtained
by substituting one or several
remainder terms $(s-t)^jr(t,s)$ into $b_j(x,\sbull)$;
using that $b_j(x,\sbull)$ is symmetric,
such summands can be written in the form
\[
h(t,s)=b_j(\eta(t),r(t,s), z_1(t,s),\ldots,z_{j-1}(t,s))\,,
\]
where each $z_i(t,s)$ is either $(s-t)^kr(t,s)$,
or $(s-t)^\ell c_\ell(t)$ for some $\ell$.
Again, any such $h\!: I\times I\to F$
is continuous, and vanishes on the diagonal.
Finally, $\rho$ involves a summand $P\!: I\times I\to F$
defined via
\[
P(t,s):=\left\{
\begin{array}{cl}
(s-t)^{-k}R(\eta(t),\eta(s)) & \mbox{if $s\not=t$;}\\
0 & \mbox{if $s=t$.}
\end{array}
\right.
\]
Then, by definition, $P$ vanishes on the diagonal.
To see that $P$ is continuous, we only need to show
that $P(s,t)\to P(s_0,s_0)=0$
whenever $(s,t)\to (s_0,s_0)$, with $s\not=t$
(cf.\ \cite{Eng}, Exerc.\,3.2.B).
But this is the case:
If $s\not=t$ are such that $\eta(t)=\eta(s)$,
then $P(t,s)=0$. If $\eta(t)\not=\eta(s)$ on the other
hand, then
\[
\|P(s,t)\|=\frac{\|R(\eta(t),\eta(s)\|}{\|\eta(t)-\eta(s)\|^k}
\cdot \left\|\frac{\eta(t)-\eta(s)}{t-s}\right\|^k\,,
\]
where the first term tends to~$0$
since $R$ is a $k$-th order remainder,
and where the second term can be written as
$\|\eta^{<1>}(t,s)\|^k$
(with notation as in Definition~\ref{defncek} above),
where $\eta^{<1>}\!: I\times I\to F$
is continuous. Hence also the product
tends to~$0$. Being a sum of continuous maps
vanishing on the diagonal, also
$\rho\!: I\times I\to F$ is continuous and
vanishes on the diagonal.
Thus (\ref{vorlage}) shows
that the curve $f\circ \eta$ admits a Taylor
expansion as described
in \cite{Sh2}, Thm.\,83.5
(resp., Prop.\,27.2\,($\gamma$),
resp., Prop.\,28.4 if $k\leq 2$),
whence $f\circ \eta$ is of class $C^k$
by the cited theorem (resp., proposition).\footnote{The full proof
of \cite{Sh2}, Thm.\,83.5
is given in
\cite{Sh1}. Only $\K$-valued
functions are considered
there, but the cited results remain valid,
with identical proofs, for maps
with values in ultrametric Banach spaces.}

Choosing $\eta$ as the curve $\eta(t):=x+ty$
for given elements
$x\in U$ and $y\in E$, we obtain
$f(x+ty)-f(x)=f(\eta(t))-f(\eta(0))=
\sum_{j=1}^k t^j a_j(x,y)+ t^k\rho(0,t)$.
Since $\rho(0,\sbull)$ is continuous and $\rho(0,0)=0$,
the coefficients $a_j(x,y)$ in the preceding formula
are uniquely determined by the function
$t\mto f(x+ty)-f(x)$, and hence by $f$ (see \cite[La.\,5.2]{BGN}).\Punkt
\noindent
{\footnotesize
{\bf Helge Gl\"{o}ckner}, TU~Darmstadt, FB~Mathematik~AG~5,
Schlossgartenstr.\,7, 64289 Darmstadt, Germany.\\
E-Mail: gloeckner@mathematik.tu-darmstadt.de}
\end{document}